\documentclass[11pt,reqno]{amsart}
\usepackage{amssymb}
\usepackage{amsmath, amssymb, amsthm}
\usepackage{mathrsfs}

\newtheorem{theorem}{Theorem}[section]
\newtheorem{definition}{Definition}[section]
\newtheorem{lemma}{Lemma}[section]
\newtheorem{remark}{Remark}[section]

\newtheorem{corollary}{Corollary}[section]
\numberwithin{equation}{section}

\newcommand{\supp}{{\mathrm {supp}}}

\begin{document}
\title[Navier-Stokes equations]{Existence results for viscous polytropic fluids with degenerate viscosity coefficients and vacuum}

\author{Shengguo Zhu}
\address[S. G.  Zhu]{Department of Mathematics, Shanghai Jiao Tong University,
Shanghai 200240,  P.R.China.}
\email{\tt zhushengguo@sjtu.edu.cn}

\begin{abstract}
In this paper, we considered the  isentropic Navier-Stokes equations for compressible fluids with density-dependent viscosities in $\mathbb{R}^3$. These systems come  from the Boltzmann equations through the Chapman-Enskog expansion to the second order, cf.\cite{tlt}, and  are degenerate when vacuum appears. We firstly  establish the  existence of the   unique local regular solution  (see Definition \ref{d1} or \cite{sz3}) when the initial data are arbitrarily large with vacuum at least appearing in the far field.  Moreover  it is interesting to show that we could't obtain any global regular solution that the $L^\infty$ norm of $u$ decays to zero as time $t$ goes to infinity. 
\end{abstract}

\date{Aug.21, 2014}

\maketitle

\section{Introduction}\ \\

 Our model is motivated by the physical consideration that in the derivation of the Navier-Stokes equations from the Boltzmann equations through the Chapman-Enskog expansion to the second order, cf.\cite{tlt}, the viscosities are not constants but depend on temperature. In particular, the viscosities of gas are proportional to the square root of the temperature for hard sphere collision. For isentropic flow, this dependence is reduced into the dependence on  density by the laws of Boyle and Gay-Lussac for ideal gas.
So the compressible isentropic Navier-Stokes equations (CINS) with degenerate viscosities in  $\mathbb{R}^3$  can be written as
\begin{equation}
\label{eq:1.1}
\begin{cases}
\rho_t+\text{div}(\rho u)=0,\\[4pt]
(\rho u)_t+\text{div}(\rho u\otimes u)
  +\nabla
   P =\text{div} \mathbb{T}.
\end{cases}
\end{equation}
We look for local strong solution with initial data
\begin{equation} \label{eq:2.211m}
\begin{split}
&(\rho,u)|_{t=0}=(\rho_0(x),u_0(x)),\quad x\in \mathbb{R}^3,
%\ \text{in} \ \mathbb{R}^3, \ (H,u)=0 \ \text{on}\ (0,T)\times \partial\mathbb{R}^3,\\
%&(H(t,x), u(t,x))\rightarrow 0, \ \rho(t,x)\rightarrow \rho^\infty, \ \text{as} \ |x|\rightarrow \infty,\ \text{for}\ (t,x)\in (0,T)\times \mathbb{R}^3.
\end{split}
\end{equation}
and far field behavior
\begin{equation} \label{eq:2.211}
\begin{split}
(\rho,u)\rightarrow (0,0) \quad \text{as } \quad |x|\rightarrow \infty,\quad t> 0.
\end{split}
\end{equation}

In system (\ref{eq:1.1}), $x\in \mathbb{R}^3$ is the spatial coordinate; $t\geq 0$ is the time;  $\rho$ is the density; $u=(u^{(1)},u^{(2)}, u^{(3)})^\top\in \mathbb{R}^3$ is the velocity of fluids; we only study the polytropic fluid, so the pressure  $P$ has the following form
\begin{equation}
\label{eq:1.2}
P=A\rho^{\gamma}, \quad 1<\gamma\leq 3,
\end{equation}
where $A$ is a positive constant, $\gamma$ is the adiabatic index. $\mathbb{T}$ is the stress tensor given by
\begin{equation}
\label{eq:1.3}
\mathbb{T}=\mu(\rho)(\nabla u+(\nabla u)^\top)+\lambda(\rho) \text{div}u\mathbb{I}_3,
\end{equation}
 where $\mathbb{I}_3$ is the $3\times 3$ unit matrix, $\mu(\rho)=\alpha\rho$ is the shear viscosity, $\lambda(\rho)=\rho E(\rho)$ is the second viscosity, where the constant $\alpha$ and function $E(\rho)$  satisfy
 \begin{equation}\label{10000}\alpha>0,\quad  2\alpha+3 E(\rho)\geq 0, \quad \text{and}\quad E(\rho)\in C^2(\overline{\mathbb{R}}^+).
 \end{equation}
For example, we can choose $\mu=\rho$ and $\lambda(\rho)=\rho^b$ for  $b= 1, 2$ or any $b\geq 3$.

%In the theory of gas dynamics, considering  the derivation of the Navier-Stokes equations from the Boltzmann equations through the Chapman-Enskog expansion to the second order, cf. \cite{bd} and \cite{tlt}, $(\mu,\lambda)$ are shown to be functions of absolute temperature in a 
%power law. If we restrict the gas flow to be isentropic, such dependence is inherited through the laws of Boyle and Gay-Lussac, and one finds that both  viscosities  are proportional to powers of density, see \cite{tlt}.  

When the  initial density has positive lower bound, the local existence of classical solutions for \eqref{eq:1.1}--\eqref{eq:2.211m} follows from a standard Banach fixed point argument due to the contraction property of the solution operators of the linearized problem, c.f. \cite{nash}.  However, when the density function connects to vacuum continuously,  this approach is not applicable for our system (\ref{eq:1.1}) due to the degeneracies caused by vacuum.  Generally it cannot be avoided    when some physical requirements are imposed, such as finite total  mass and  energy in the whole space $\mathbb{R}^3$, because at least we need that 
$$\rho(t,x)\rightarrow 0,\quad \text{as}\quad |x|\rightarrow +\infty.$$

When  $(\mu,\lambda)$ are both constants, for the existence of  $3$D solutions of the  isentropic flow with  arbitrary data, the main breakthrough is due to Lions \cite{lions}, where he established the global existence of weak solutions  in $\mathbb{R}^3$, periodic domains or bounded domains with homogenous Dirichlet boundary conditions provided $\gamma >9/5$. The restriction on $\gamma$ is improved to $\gamma>3/2$ by Feireisl \cite{fu1}\cite{fu2}, and the corresponding result for the non-isentropic flow can be seen in \cite{fu3}.  Recently in Cho-Choe-Kim \cite{CK3}\cite{CK}, via introducing the following initial layer compatibility condition:
\begin{equation*}
-\text{div} \mathbb{T}_{0}+\nabla P(\rho_0)=\sqrt{\rho_{0}} g
\end{equation*}
for some $g\in L^2$,   a local theory for arbitrarily large strong solutions was established successfully; see also \cite{luoluo}.
And Huang-Li-Xin \cite{HX1} obtained the global well-posedness of classical solutions with small energy and vacuum to Cauchy problem  for isentropic flow.

When $(\mu, \lambda)$  are both dependent of $\rho$ as shown in the following form:
\begin{equation}
\label{fandan}
\mu(\rho)=\alpha\rho^{\delta_1},\quad \lambda(\rho)=\beta\rho^{\delta_2},
\end{equation}
where  $\delta_1>0$, $\delta_2>0$, $\alpha>0$ and $\beta$ are all real constants, system (\ref{eq:1.1}) has received a lot of attention recently, see \cite{bd2}\cite{bd}\cite{bd3}\cite{hailiang}\cite{taiping}\cite{vassu}\cite{tyc2}\cite{zyj}. However,  except for the $1$D problems, there are still only  few results   on  the   strong solutions  for the multi-dimensional problems because of the possible degenerency for the Lam$\acute{\text{e}}$  operator  caused by the initial vacuum. This degeneracy gives rise to some difficulties in the regularity estimate because of the less regularizing effect of the viscosity on  solutions.  This is one of the major obstacles preventing us from utilizing a similar remedy proposed by Cho et. al. for the case of constant viscosity coefficients. However, recently in $2$D  space, Li-Pan-Zhu \cite{sz3} has obtained the existence of the unique local classical solutions for system (\ref{eq:1.1}) under the assumptions 
$$\rho_0\rightarrow 0,\quad  \text{as}\quad  |x|\rightarrow \infty$$
 and 
\begin{equation}\label{100000}\delta_1=1, \quad \delta_2=0 \ \text{or}\  1,\quad \alpha>0,\quad \alpha+\beta\geq 0,
\end{equation}
 but the vacuum cannot appear in any local point. And in \cite{sz345}, they also  proved the existence of the unique local classical solutions for system (\ref{eq:1.1}) under the assumption
$$1<\delta_1=\delta_2<\min \Big(3,\frac{\gamma+1}{2}\Big), \quad \alpha>0,\quad \alpha+\beta\geq 0$$
 with  initial vacuum  appearing in some open set or the far field. 

In this paper,  we  generalize the $2$D existence  result obtained in \cite{sz3} to $\mathbb{R}^3$ in $H^2$ space and    assume  (\ref{10000}) instead of (\ref{fandan})-(\ref{100000}).  Moreover,  we will show an very interestring phenomenon  that that it is impossible to  obtain any global regular solution that the $L^\infty$ norm of $u$ decays to zero as time $t$ goes to infinity. 

%,  and even the short time well-posedness of classical solutions with $\inf \rho_0=0$ opens for $0<\delta_1 <1$

Throughout this paper, we adopt the following simplified notations for the standard homogeneous and inhomogeneous Sobolev space:
\begin{equation*}\begin{split}
&D^{k,r}=\{f\in L^1_{loc}(\mathbb{R}^3): |f|_{D^{k,r}}=|\nabla^kf|_{L^r}<+\infty\},\quad D^k=D^{k,2}(k\geq 2), \\[6pt]
&D^{1}=\{f\in L^6(\mathbb{R}^3): |f|_{D^{1}}=|\nabla f|_{L^2}<\infty\},\quad \|(f,g)\|_X=\|f\|_X+\|g\|_X,\\[6pt]
 &\|f\|_s=\|f\|_{H^s(\mathbb{R}^3)},\quad |f|_p=\|f\|_{L^p(\mathbb{R}^3)},\quad |f|_{D^k}=\|f\|_{D^k(\mathbb{R}^3)}.
\end{split}
\end{equation*}
A detailed study of homogeneous Sobolev space can be found in \cite{gandi}.

First we introduce the definitions of regular solutions and strong solutions to Cauchy problem (\ref{eq:1.1})-(\ref{eq:2.211}). Via  introducing the new variable $c(t,x)=\sqrt{A\gamma}\rho^{\frac{\gamma-1}{2}}$   (local sound speed) and $\psi=\frac{2}{\gamma-1}\nabla c / c=(\psi^{(1)},\psi^{(2)},\psi^{(3)})^\top$, then  (\ref{eq:1.1})-(\ref{eq:2.211})  can be written as
\begin{equation}
\begin{cases}
\label{eq:cccq}
\displaystyle
c_t+u\cdot \nabla c+\frac{\gamma-1}{2}c\text{div} u=0,\\[8pt]
\displaystyle
%\psi_t+\nabla (u\cdot \psi)+\nabla \text{div} u=0,\\[8pt]
\displaystyle
u_t+u\cdot\nabla u +\frac{2}{\gamma-1}c\nabla c+Lu=\psi\cdot Q(c,u),\\[8pt]
(c,u)|_{t=0}=(c_0,u_0),\quad x\in \mathbb{R}^3,\\[8pt]
(c,u)\rightarrow (0,0) \quad \text{as } \quad |x|\rightarrow \infty,\quad t> 0,
 \end{cases}
\end{equation}
where $L$ is the so-called Lam$\acute{ \text{e} }$ operator given by
$$ Lu=-\text{div}(\alpha(\nabla u+(\nabla u)^\top)+ \overline{E}(c)\text{div}u\mathbb{I}_3),$$
and terms $\big( Q(c,u), \overline{E}(c)\big)$ are given by
$$ Q(c,u)=\alpha(\nabla u+(\nabla u)^\top)+ \overline{E}(c)\text{div}u\mathbb{I}_3,\quad \overline{E}(c)=E\big(((A\gamma)^{\frac{-1}{2}}c)^{\frac{2}{\gamma-1}}\big).$$
Similar to \cite{sz3},  the regular solution is defined via:
\begin{definition}[\text{\textbf{Regular solutions to Cauchy problem (\ref{eq:1.1})-(\ref{eq:2.211})}}]\label{d1}\ \\[2pt]
 Let $T> 0$ be a finite constant. $(c,u)$  is called a regular solution to Cauchy problem (\ref{eq:1.1})-(\ref{eq:2.211})  in $ [0,T]\times \mathbb{R}^3$ if $(c,u)$ satisfies
\begin{equation*}\begin{split}
&(\textrm{A})\quad (c,u)\   \text{satisfies the Cauchy problem (\ref{eq:cccq})  a.e. in  }\ (t,x)\in (0,T]\times \mathbb{R}^3; \\
&(\textrm{B})\quad c \geq 0,\quad  c\in C([0,T]; H^2), \ c_t \in C([0,T]; H^1); \\
&(\textrm{C})\quad  \psi\in C([0,T] ; D^1),\ \psi_t\in C([0,T] ; L^2);\\
&(\textrm{D})\quad u\in C([0,T]; H^2)\cap L^2([0,T] ;  D^3), \ u_t \in C([0,T]; L^2)\cap L^2([0,T] ; D^1).
\end{split}
\end{equation*}
%\begin{equation} \label{eq:5.1}\end{equation}
\end{definition}

This definition for regular solutions is similar to that of  Makino-Ukai-Kawashima \cite{tms1}, which studied the local existence of classical solutions to non-isentropic Euler equations with initial data arbitrarily large and $\inf \rho_0=0$. Some similar definitions can also be seen in \cite{sz3}\cite{sz345}\cite{sz2}\cite{dcds}\cite{tpy}\cite{tms1}\cite{makio}\cite{tyc2}. And the strong solution can be given as

\begin{definition}[\text{\textbf{Strong solutions to Cauchy problem (\ref{eq:1.1})-(\ref{eq:2.211})}}]\label{d2}\ \\[2pt]
 Let $T> 0$ be a finite constant. $(\rho,u)$  is called a strong solution to Cauchy problem (\ref{eq:1.1})-(\ref{eq:2.211})  in $ [0,T]\times \mathbb{R}^3$ if $(\rho,u)$ satisfies
\begin{equation*}\begin{split}
&(\textrm{A1})\quad (\rho,u)\   \text{satisfies the Cauchy problem (\ref{eq:1.1})-(\ref{eq:2.211})  a.e. in  }\ (t,x)\in (0,T]\times \mathbb{R}^3; \\
&(\textrm{B1})\quad \rho\geq 0,\ \rho\in C([0,T]; H^2), \ \rho_t \in C([0,T]; H^1); \\
&(\textrm{C1})\quad u\in C([0,T]; H^2)\cap L^2([0,T] ; D^3), \ u_t \in C([0,T]; L^2)\cap L^2([0,T] ; D^1);\\
&(\textrm{D1})\quad u_t+u\cdot\nabla u +Lu=(\nabla \rho/\rho) \cdot Q(c,u) \quad  \text{holds when }\ \rho(t,x)=0.
\end{split}
\end{equation*}
%\begin{equation} \label{eq:5.1}\end{equation}
\end{definition}
\begin{remark}\label{zhen90}
It is obvious that   conditions $(\textrm{B})$ or $(\textrm{B1})$   mean  that the vacuum must  appear at least  in the far field. 

\end{remark}

Now we give the main existence results of this paper:
\begin{theorem}[\textbf{Existence of the unique local regular solution}]\label{th2}\ \\[2pt]
Let $1< \gamma \leq 2$ or $\gamma = 3$. If the initial data $( c_0, u_0)$ satisfies the regularity condition
\begin{equation}\label{th78}
\begin{split}
&c_0\geq 0,\quad (c_0, u_0)\in H^2, \quad \psi_0\in  D^1,
\end{split}
\end{equation}
 then there exists a small time $T_*$ and a unique regular solution $(c, u)$ to Cauchy problem (\ref{eq:1.1})-(\ref{eq:2.211}).
Moreover,  we also have $\rho(t,x)\in C([0,T_*]\times \mathbb{R}^3)$.
\end{theorem}

\begin{remark}\label{r2} First we remark that \eqref{th78} identifies a class of admissible initial data that 
provides unique solvability to our problem \eqref{eq:1.1}--\eqref{eq:2.211}. On the other hand, this set 
of initial data contains a large class of functions, for example,
$$
\rho_0(x)=\frac{1}{1+|x|^{2\sigma}}, \quad u_0(x)=0,\quad x\in \mathbb{R}^3,
$$
where $\sigma>\max\{1, \frac{1}{\gamma-1}\}$. 

Second, we remark that under the initial assumption (\ref{th78}) and $\rho^{b-1}_0\in H^2$, the conclusion obtained in  Theorem \ref{th2} still holds for the case that $\lambda(\rho)=\rho^b$ (i.e., $E(\rho)=\rho^{b-1}$) when $b\in (1,2) \cup (2,3)$ and $1<\gamma \leq 3$. The details can be seen in Subsection 3.5.
\end{remark}

According to the conclusions obtained in Theorem \ref{th2} and the standard quasi-linear hyperbolic equations theory, we quickly have the following result:
\begin{corollary}[\textbf{Existence of strong solutions}]\label{co2}\ \\[2pt]
Let $1< \gamma \leq 2$ or $\gamma = 3$. Then the regular solution obtained in Theorem \ref{th2} is indeed the strong solution to Cauchy problem (\ref{eq:1.1})-(\ref{eq:2.211}).
\end{corollary}

Next, we will show some interesting phenomenon which  tells us that there does not exist any global regular solution to Cauchy problem  (\ref{eq:1.1})-(\ref{eq:2.211}) with the $L^\infty$  norm of  velocity $u$  decaying to  zero as time goes to infinity. Let
\begin{align*}
&\mathbb{P}(t)=\int_{\mathbb{R}^{3}}\rho(t,x)u(t,x)\text{d}x \quad \textrm{(total momentum)}.
\end{align*}
\begin{theorem}[\textbf{Non-existence of global solutions with $L^\infty$ decay on $u$}]\label{th:2.20}\ \\[2pt]
 Let $1< \gamma \leq 2$. Adding   $0<|\mathbb{P}(0)|$ to (\ref{th78}).
Then there is no global  regular  solution $(\rho,u)$ obtained in Theorems  \ref{th2}.
satisfying the following decay
\begin{equation}\label{eq:2.15}
\limsup_{t\rightarrow +\infty} |u(t,x)|_{\infty}=0.
\end{equation}

\end {theorem}

However,  via combining  the arguments used in this paper and  \cite{sz3}   in $\mathbb{R}^2$, we can also have the similar conclusions  obtained above in $H^2$ space:
\begin{theorem}\label{th22} Let $1< \gamma \leq 2$ or $\gamma = 3$. If the initial data $( \rho_0, u_0)$ satisfy 
$$
0\leq \rho^{\frac{\gamma-1}{2}}_0\in H^2(\mathbb{R}^2),\quad u_0\in H^2(\mathbb{R}^2),\quad \nabla \rho_0/\rho_0 \in L^6(\mathbb{R}^2)\cap D^1(\mathbb{R}^2),
$$
 then there exists a time $T_*>0$ and a unique regular solution $(\rho, u)$ to the Cauchy problem (\ref{eq:1.1})-(\ref{eq:2.211}) satisfying
\begin{equation}\label{reg11}\begin{split}
& \rho^{\frac{\gamma-1}{2}} \in C([0,T_*];H^2(\mathbb{R}^2)),\ (\rho^{\frac{\gamma-1}{2}})_t \in C([0,T_*];H^1(\mathbb{R}^2)),\\
& \nabla \rho/\rho \in C([0,T_*]; L^6\cap D^1(\mathbb{R}^2)),\quad  (\nabla \rho/\rho)_t \in C([0,T_*]; L^2(\mathbb{R}^2)),\\
& u\in C([0,T_*]; H^2(\mathbb{R}^2))\cap L^2([0,T_*] ; D^{3}(\mathbb{R}^2)),\\
& u_t \in C([0,T_*]; L^2(\mathbb{R}^2))\cap L^2([0,T_*] ; D^1(\mathbb{R}^2)).
%\ t^{\frac{1}{2}}u \in L^\infty([0,T_*]; D^3(\mathbb{R}^2)).
%& t^{\frac{1}{2}}u_t \in L^\infty([0,T_*]; D^1(\mathbb{R}^2))\cap L^2([0,T_*];D^2(\mathbb{R}^2)), \ t^{\frac{1}{2}}u_{tt}\in L^2([0,T_*];L^2(\mathbb{R}^2)).
\end{split}
\end{equation}
Moreover, we also have 
$
\rho(t,x)\in C([0,T_*]\times \mathbb{R}^3)
$,
and
$$
 \rho\in C([0,T_*];H^2(\mathbb{R}^2)),\quad  \rho_t \in C([0,T_*];H^1(\mathbb{R}^2)).
$$
\end{theorem}

The rest of  this paper is organized as follows. In Section $2$, we  give some important lemmas that will be used frequently in our proof.
In Section $3$, we prove the  existence of the unique regular solution shown in Theorem \ref{th2}  via establishing  some a prior estimates which are independent of the lower bound of $c$, and these estimates can be  obtained by  the approximation process from non-vacuum to vacuum. In Section $4$, based on the conclusions obtained  in  Section $3$,  we  give the proof for our main result: the local existence of strong solutions to the original problem
(\ref{eq:1.1})-(\ref{eq:2.211}) shown in Corollary \ref{co2}.  Finally,  in Section $5$, we will show the non-existence of global solutions with $L^\infty$ decay on $u$.

\section{Preliminary}
In this section, we show some important lemmas  that will be frequently used in our proof.
The first one is the  well-known Gagliardo-Nirenberg inequality.
\begin{lemma}\cite{oar}\label{lem2as}\
For $p\in [2,6]$, $q\in (1,\infty)$, and $r\in (3,\infty)$, there exists some generic constant $C> 0$ that may depend on $q$ and $r$ such that for 
$$f\in H^1(\mathbb{R}^3),\quad \text{and} \quad  g\in L^q(\mathbb{R}^3)\cap D^{1,r}(\mathbb{R}^3),$$
 we have
\begin{equation}\label{33}
\begin{split}
&|f|^p_p \leq C |f|^{(6-p)/2}_2 |\nabla f|^{(3p-6)/2}_2,\\[8pt]
&|g|_\infty\leq C |g|^{q(r-3)/(3r+q(r-3))}_q |\nabla g|^{3r/(3r+q(r-3))}_r.
\end{split}
\end{equation}
\end{lemma}
Some common versions of this inequality can be written as
\begin{equation}\label{ine}\begin{split}
|u|_6\leq C|u|_{D^1},\quad |u|_{\infty}\leq C|u|^{\frac{1}{2}}_6|\nabla u|^{\frac{1}{2}}_6\leq  C(|u|_{D^1}+|u|_{D^2}), \quad |u|_{\infty}\leq C\|u\|_{W^{1,r}}.
\end{split}
\end{equation}

The second one  can be seen in Majda \cite{amj}, here we omit its proof.

\begin{lemma}\cite{amj}\label{zhen1}
Let constants $r$, $a$ and $b$ satisfy the relation 
$$\frac{1}{r}=\frac{1}{a}+\frac{1}{b},\quad \text{and} \quad 1\leq a,\ b, \ r\leq \infty.$$  $ \forall s\geq 1$, if $f, g \in W^{s,a} \cap  W^{s,b}(\mathbb{R}^3)$, then we have
\begin{equation}\begin{split}\label{ku11}
&|\nabla^s(fg)-f \nabla^s g|_r\leq C_s\big(|\nabla f|_a |\nabla^{s-1}g|_b+|\nabla^s f|_b|g|_a\big),
\end{split}
\end{equation}
\begin{equation}\begin{split}\label{ku22}
&|\nabla^s(fg)-f \nabla^s g|_r\leq C_s\big(|\nabla f|_a |\nabla^{s-1}g|_b+|\nabla^s f|_a|g|_b\big),
\end{split}
\end{equation}
where $C_s> 0$ is a constant only depending on $s$.
\end{lemma}

Based on  harmonic analysis,  we introduce a regularity estimate result for the following elliptic problem in the whole domain $\mathbb{R}^3$:
\begin{equation}\label{ok}
-\triangle u=f, \quad  u\rightarrow 0 \quad \text{as} \ |x|\rightarrow \infty.
\end{equation}
\begin{lemma}\cite{harmo}\label{zhenok}
If $u\in D^{1,p}$ with $1< p< \infty$ is a weak solution to system (\ref{ok}), then
$$
\|u\|_{D^{2,p}(\mathbb{R}^3)} \leq C \|f\|_{L^p(\mathbb{R}^3)},
$$
where $C$ depending only on  $p$. Moreover, if $f=\text{div} \  h$, then we also have
$$
\|u\|_{D^{1,p}(\mathbb{R}^3)} \leq C \|h\|_{L^p(\mathbb{R}^3)}.
$$
\end{lemma}
\begin{proof}
The  proof can be obtained via the classical  harmonic analysis \cite{harmo}.
\end{proof}

Finally, the  last one is some result obtained via the Aubin-Lions Lemma.
\begin{lemma}\cite{jm}\label{aubin} Let $X_0$, $X$ and $X_1$ be three Banach spaces with $X_0\subset X\subset X_1$. Suppose that $X_0$ is compactly embedded in $X$ and that $X$ is continuously embedded in $X_1$.\\[1pt]

I) Let $G$ be bounded in $L^p(0,T;X_0)$ where $1\leq p < \infty$, and $\frac{\partial G}{\partial t}$ be bounded in $L^1(0,T;X_1)$. Then $G$ is relatively compact in $L^p(0,T;X)$.\\[1pt]

II) Let $F$ be bounded in $L^\infty(0,T;X_0)$  and $\frac{\partial F}{\partial t}$ be bounded in $L^p(0,T;X_1)$ with $p>1$. Then $F$ is relatively compact in $C(0,T;X)$.
\end{lemma}

\section{Existence of the unique regular solutions}

In this section, we will give the proof for the existence of the unique regular solutions shown in Theorem \ref{th2} by Sections $3.1$-$3.4$.
\subsection{Linearization}
 For simplicity, in the following sections, we denote $\frac{1}{\gamma-1}=\theta$.
Now we consider the following linearized equations
\begin{equation}\label{li4}
\begin{cases}
c_t+v\cdot \nabla c+\frac{\gamma-1}{2}c\text{div} v=0,\\[8pt]
%\psi_t+\nabla (v\cdot \psi)+\nabla \text{div} v=0,\\[8pt]
u_t+v\cdot\nabla v +2\theta c\nabla c+Lu=\psi\cdot Q(c,v),
 \end{cases}
\end{equation}
 where $\psi=2\theta\nabla c/c$ and
\begin{equation} \label{eq:5.231}  Q(c,v)=\alpha(\nabla v+(\nabla v)^\top)+ \overline{E}(c)\text{div}v\mathbb{I}_3. \end{equation}
The initial data is given by
\begin{equation}\label{qwe1}
(c,\psi,u)|_{t=0}=(c_0,\psi_0,u_0),\quad x\in \mathbb{R}^3.
\end{equation}
We assume that
\begin{equation}\label{th78rr}
\begin{split}
&c_0\geq 0,\quad  (c_0-c^\infty, u_0)\in H^2, \quad  \psi_0=2\theta\nabla c_0/c_0\in D^1
\end{split}
\end{equation}
where $c^\infty\geq 0$ is a constant.  And $v=(v^{(1)},v^{(2)},v^{(3)})^\top\in \mathbb{R}^3$ is a known vector satisfying
\begin{equation}\label{vg}
\begin{split}
& v\in C([0,T] ; H^2)\cap L^2([0,T] ; D^{3}), \  v_t \in C([0,T] ; L^2)\cap L^2([0,T] ; D^1).
%\ t^{\frac{1}{2}}v \in L^\infty([0,T_*]; D^3).
%&t^{\frac{1}{2}}v_t \in L^\infty([0,T_*]; D^1)\cap L^2([0,T_*];D^2),\ t^{\frac{1}{2}}v_{tt}\in L^2([0,T];L^2).
\end{split}
\end{equation}
Moreover, we assume that $u_0=v(t=0,x)$.
Then we have  the following existence of a strong solution $(c, \psi ,u)$ to (\ref{li4})-(\ref{vg}) by the standard methods at least in the case that the initial data is away from vacuum.

 \begin{lemma}\label{lem1}
 Assume that the initial data (\ref{qwe1}) satisfy (\ref{th78rr}) and $c_0> \delta$ for some positive constant.
Then there exists a unique strong solution $(c, \psi ,u)$ to (\ref{li4})-(\ref{vg}) such that
\begin{equation}\label{reggh}\begin{split}
&c\geq \underline{\delta}, \ c-c^\infty \in C([0,T]; H^2), \ c _t \in C([0,T]; H^1), \\
&\psi \in C([0,T] ; D^1),\ \psi_t \in C([0,T] ; L^2) ,\\
& u\in C([0,T]; H^2)\cap L^2([0,T]; D^{3}), \ u_t \in C([0,T]; L^2)\cap L^2([0,T] ; D^1),
%\ t^{\frac{1}{2}}u\in L^\infty([0,T];D^3),\\
%&  t^{\frac{1}{2}}u_t \in L^\infty([0,T_*]; D^1)\cap L^2([0,T_*];L^2),\ t^{\frac{1}{2}}u_{tt}\in L^2([0,T_*]; L^2).
\end{split}
\end{equation}
where $\underline{\delta}$ is a positive constant.
\end{lemma}
\begin{proof}
First, the existence  of the solution $c$ to $(\ref{li4})_1$ can be obtained essentially via Lemma 6 in  \cite{CK} via the standard hyperbolic theory. And $c$ can be written as
\begin{equation}
\label{eq:bb1}
c(t,x)=c_0(U(0;t,x))\exp\Big(-\frac{\gamma-1}{2}\int_{0}^{t}\textrm{div}v(s,U(s;t,x))\text{d}s\Big),
\end{equation}
where  $U\in C([0,T]\times[0,T]\times \mathbb{R}^3)$ is the solution to the initial value problem
\begin{equation}
\label{eq:bb1}
\begin{cases}
\frac{d}{dt}U(t;s,x)=v(t,U(t;s,x)),\quad 0\leq t\leq T,\\[6pt]
U(s;s,x)=x, \quad \ \quad \quad 0\leq s\leq T,\ x\in \mathbb{R}^3.
\end{cases}
\end{equation}
So we easily know that there exists a positive constant $\underline{\delta}$ such that $c\geq \underline{\delta}$.

Second, due to  $c\geq \underline{\delta}$, we quickly obatin that 
$$
\psi \in C([0,T] ; D^1),\ \psi_t \in C([0,T] ; L^2) .
$$

At last, based on the regularity of $c$ and $\psi$, the desired conclusions for $u$ can be obtained from the linear parabolic equations
$$
u_t+v\cdot\nabla v +2\theta c\nabla c+Lu=\psi\cdot Q(c,v)
$$
via the classical Galerkin methods which can be seen in \cite{CK3}\cite{CK}, here we omit it.
\end{proof}

\subsection{A prior estimate}  In this section, we assume that $(c,\psi,u)$ is the unique strong solution to (\ref{li4})-(\ref{vg}), then we will get some   a prior estimates which are independent of the lower bound $\delta$ of $c_0$.
Now we fix a positive constant $c_0$ large enough that
\begin{equation}\label{houmian}\begin{split}
2+c^\infty+|c_0|_{\infty}+\|c_0-c^\infty\|_{2}+|\psi_0|_{D^1}+\|u_0\|_{2}\leq b_0,
\end{split}
\end{equation}
and
\begin{equation}\label{jizhu1}
\begin{split}
\sup_{0\leq t \leq T^*}| v(t)|^2_{2}+\int_{0}^{T^*}  |\nabla v(t)|^2_{2} \text{d}t \leq& b^2_1,\\
\sup_{0\leq t \leq T^*}| v(t)|^2_{D^1}+\int_{0}^{T^*} \Big( | v(t)|^2_{D^2}+|v_t(t)|^2_{2}\Big)\text{d}t \leq& b^2_2,\\
\sup_{0\leq t \leq T^*}\big(|v(t)|^2_{D^2}+|v_t(t)|^2_{2}\big)+\int_{0}^{T^*} \Big( |v(t)|^2_{D^3}+|v_t(t)|^2_{D^1}\Big)\text{d}t \leq& b^2_3
\end{split}
\end{equation}
for some time $T^*\in (0,T)$ and constants $b_i$ ($i=1,2,3$) such that 
$1< b_0\leq b_1 \leq b_2\leq b_3 $. The constants $b_i$ ($i=1,2,3$)  and time $T^*$ will be determined later and depend only on $b_0$,  the fixed constants $\alpha$, $A$, $\gamma$ and  $T$ (see (\ref{dingyi})).
Throughout this and next two subsections, we denote by $C$ a generic positive constant depending only on fixed constants $\alpha$, $A$, $\gamma$ and $T$.
Moreover, let $1\leq M(\cdot)\in C(\overline{\mathbb{R}}^+)$ be  a nondecreasing and continuous function, which only dependes on $E(\cdot)$ and the  constant $C$.
To begin with, we give some estimates for $c$.

\begin{lemma}[\textbf{Estimates for $c$}]\label{2}
$$|c(t)|^2_\infty+\|c(t)-c^\infty\|^2_{2}\leq Cb^2_0,\quad |c_t(t)|_{2}\leq Cb_0b_2,\quad |c_t|_{D^1}\leq Cb_0b_3,$$
$$|\overline{E}(c)(t)|^2_\infty+\|\overline{E}(c)(t)-\overline{E}(c^\infty)\|^2_{2}\leq M(b_0),$$
$$|\overline{E}(c)_t(t)|_{2}\leq M(b_0)b_0b_2,\quad |\overline{E}(c)_t(t)|_{D^1}\leq M(b_0)b_0b_3$$
for $0\leq t \leq T_1=\min (T^*, (1+b_3)^{-2})$.
\end{lemma}

\begin{proof}\underline{Step 1}. 
From stand energy estimate theories introduced in \cite{CK},  we easily have
\begin{equation*}
\begin{split}
\|c(t)-c^\infty\|_{2}\leq \Big(\|c_0-c^\infty\|_{2} +c^\infty\int_0^t \|\nabla v(s)\|_{2}\text{d}s\Big)\exp\Big(C\int_0^t    \| \nabla v(s)\|_{ 2 } \text{d}s\Big).
\end{split}
\end{equation*}
Therefore, observing that
$$
\int_0^t \|\nabla v(s)\|_{2}\text{d}s\leq t^{\frac{1}{2}}\Big(\int_0^t \| \nabla v(s)\|^2_{2}\text{d}s\Big)^{\frac{1}{2}}\leq C(b_2t+b_3t^{\frac{1}{2}}),
$$
then the estimate for $\|c-c^{\infty}\|_{2}$ is available for $0\leq t \leq T_1=\min (T^*, (1+b_3)^{-2})$.

The estimate for $c_t$   follows  from the following relation
$$c_t=-v\cdot \nabla c-\frac{\gamma-1}{2}c\text{div} v,$$
we easily have,  for $0\leq t \leq T_1$,
\begin{equation}\label{zhen6}
\begin{cases}
|c_t(t)|_{2}\leq C\big(|v(t)|_6| \nabla c(t)|_{3}+|c(t)|_\infty|\text{div} v(t)|_{2}\big)\leq Cb_0b_2,\\[10pt]
|c_t(t)|_{D^1}\leq C\big(|v(t)|_{\infty} |c(t)|_{D^2}+|c(t)|_{\infty}| v(t)|_{D^2}+|\nabla c(t)|_{6}| \nabla v(t)|_{3}\big)\leq Cb_0b_3.
\end{cases}
\end{equation}

 \underline{Step 2}. Due to $1<\gamma \leq 2$ or $\gamma=3$,  and $ E(\rho)\in C^2(\overline{\mathbb{R}}^+)$, then we quickly know that 
$$\overline{E}(c)=E\big(((A\gamma)^{\frac{-1}{2}}c)^{\frac{2}{\gamma-1}}\big) \in C^2(\overline{\mathbb{R}}^+).$$
So the desired estimates for $\overline{E}(c)$ follows quickly from  the estimates on $c$.
\end{proof}

Next, we give some very important estimates for $\psi$.
Due to 
$$\psi=\frac{2}{\gamma-1}\nabla \phi/\phi,\quad \text{and} \quad \phi\geq \underline{\delta},$$  from $(\ref{li4})_1$ we deduce that $\psi$ satisfies
$$
\psi_t+\nabla (v\cdot \psi)+\nabla \text{div} v=0, \ \psi_0=\frac{2}{\gamma-1}\nabla \phi_0/\phi_0\in  D^1.
$$
A direct calculation shows that  
$$\partial_i \psi^{(j)}=\partial_j \psi^{(i)} \quad \text{for}\ i,j=1,2,3$$
 in distribution sense,  then the above Cauchy problem can be written as
\begin{equation}\label{ku}
\psi_t+\sum_{l=1}^3 A_l \partial_l\psi+B\psi+\nabla \text{div} v=0,\ \psi_0 \in  D^1,
\end{equation}
where 
$$A_l=(a^l_{ij})_{3\times 3},\quad \text{for}\ i,j,l=1,2,3$$
are symmetric  with 
$$a^l_{ij}=v^{(l)}\quad \text{for}\ i=j;\quad \text{otherwise}\  a^l_{ij}=0, $$
 and $B=(\nabla v)^\top$, which means that (\ref{ku}) is a positive  symmetric hyperbolic  system,
then we have the following a prior estimate for $\psi$  via the stand energy estimate theory for positive symmetric hyperbolic system. This lemma will be used to  deal with the degenerate Lam$\acute{ \text{e}}$ operator when vacuum appears for our  reformulated system.
\begin{lemma}[\textbf{Estimates for $\psi$}]\label{3}
$$|\psi(t)|^2_{D^1}\leq Cb^2_0,\quad  |\psi(t)_t|^2_{2}\leq Cb^4_3, \quad 0\leq t \leq T_1.$$

\end{lemma}
\begin{proof}

According to the proof of Lemma \ref{lem1}, we know that $\psi$ has the following regualrity
$$ \psi \in C([0,T] ; D^1),\quad \psi_t\in C([0,T] ; L^2).$$
  So, let $\varsigma=(\varsigma_1,\varsigma_2,\varsigma_3)^\top$ ($ |\varsigma|=1$ and $\varsigma_i=0,1$), differentiating (\ref{ku}) $\varsigma-$times with respect to $x$, we have
\begin{equation}\label{hyp}\begin{split}
&(D^\varsigma\psi)_t+\sum_{l=1}^3 A_l \partial_lD^\varsigma \psi+BD^\varsigma\psi+D^\varsigma\nabla \text{div} v \\
=&\big(-D^\varsigma(B\psi)+BD^\varsigma \psi\big)+\sum_{l=1}^3 \big(-D^\varsigma(A_l \partial_l \psi)+A_l \partial_lD^\varsigma \psi\big)=\Theta_1+\Theta_2.
\end{split}
\end{equation}
Multiplying (\ref{hyp}) by $2D^{\varsigma}\psi$ and  integrating over $\mathbb{R}^3$, because  $A_l$ ($l=1,2,3$) are symmetric,  we easily deduce that
\begin{equation}\label{zhenzhen}\begin{split}
\frac{d}{dt}|D^\varsigma\psi|^2_2\leq C\Big(\sum_{l=1}^{3}|\partial_{l}A_l|_\infty+|B|_\infty\Big)|D^\varsigma\psi|^2_2+C(|\Theta_1 |_2+|\Theta_2|_2+\|\nabla^2 v\|_1)|D^\varsigma\psi|_2.
\end{split}
\end{equation}
Then let $r=a=2$, $b=\infty$  when $|\varsigma|=1$ in (\ref{ku22}), we easily have
\begin{equation}\label{zhen2}
|\Theta_1|_2=|D^\varsigma(B\psi)-BD^\varsigma \psi|_2\leq C|\nabla^2 v|_3|\psi|_6;
\end{equation}
 let $r=b=2$, $a=\infty$  when $|\varsigma|=1$ in (\ref{ku22}), we easily have
\begin{equation}\label{zhen2l}
|\Theta_2|_2=|D^\varsigma(A_l \partial_l \psi)-A_l \partial_lD^\varsigma \psi|_2\leq C|\nabla v|_\infty |\nabla\psi|_2.
\end{equation}

Combining (\ref{zhenzhen})-(\ref{zhen2l}) and Lemma \ref{lem2as}, we have
\begin{equation*}
\frac{d}{dt}|\psi(t)|_{D^1}\leq C\|\nabla v \|_2|\psi(t)|_{D^1}+C\|\nabla^2 v\|_1.
\end{equation*}
According to Gronwall's inequality, we have
\begin{equation*}\begin{split}
|\psi(t)|_{D^1}\leq&  \Big(|\psi_0|_{D^1}+\int_0^t \|\nabla^2 v\|_{1} \text{d}t\Big) \exp\Big(C\int_0^t \|\nabla v\|_{2} \text{d}t\Big)
\end{split}
\end{equation*}
for $0\leq t \leq T_1$. Therefore, observing that
$$
\int_0^t \| v(s)\|_{3}\text{d}s\leq t^{\frac{1}{2}}\Big(\int_0^t \| v(s)\|^2_{3}\text{d}s\Big)^{\frac{1}{2}}\leq C(b_2t+b_3t^{\frac{1}{2}}),
$$
then desired estimate for $|\psi(t)|_{D^1}$ is available for $0\leq t \leq T_1$.

Due to the following relation
\begin{equation}\label{ghtu}\psi_t=-\nabla (v \cdot \psi)-\nabla \text{div} v,
\end{equation}
combining with the Lemma \ref{lem2as}, we easily have, for $0\leq t \leq T_1$
\begin{equation*}\begin{split}
|\psi_t(t)|_2\leq& C\big(|v|_{\infty}|\psi|_{D^1}+|\nabla v|_{3}|\psi|_{6}+|v|_{D^2}\big)(t)\leq Cb^2_3.
\end{split}
\end{equation*}

\end{proof}

Now we give the  estimates for the lower order terms of the velocity $u$.

 \begin{lemma}[\textbf{Lower order estimates of the velocity $u$}]\label{4}
\begin{equation*}
\begin{split}
|u(t)|^2_{ 2}+\int_{0}^{t}|\nabla u(s)|^2_{2}\text{d}s\leq Cb^2_0
\end{split}
\end{equation*}
for $0 \leq t \leq T_2=\min(T^*,(1+M(b_0)b^{4}_3)^{-1})$.
 \end{lemma}
\begin{proof}
 Multiplying $(\ref{li4})_2$ by $u$ and integrating  over $\mathbb{R}^3$, then we have
\begin{equation}\label{zhu1}
\begin{split}
&\frac{1}{2} \frac{d}{dt}|u|^2_2+\alpha|\nabla u|^2_2+\int_{\mathbb{R}^3}(\alpha+\overline{E}(c))|\text{div} u|^2 \text{d}x\\
=&\int_{\mathbb{R}^3} \Big(-v\cdot \nabla v \cdot u-2\theta c \nabla c \cdot u+\psi \cdot Q(c,v)\cdot u \Big) \text{d}x\equiv: \sum_{i=1}^3 I_i.
\end{split}
\end{equation}
According to H\"older's inequality, Lemma \ref{lem2as} and Young's inequality, we have
\begin{equation}\label{zhu2}
\begin{split}
I_1=&-\int_{\mathbb{R}^3} v\cdot \nabla v \cdot u \text{d}x\leq C|v|_3|\nabla v|_2|u|_6\leq C|v|^2_3|\nabla v|^2_2+\frac{\alpha}{10}|\nabla u|^2_2,\\
%\end{split}
%\end{equation}
%\begin{equation}\label{zhuoop}
%\begin{split}
I_2=&-\int_{\mathbb{R}^3} 2\theta c \nabla c\cdot u \text{d}x
\leq C|\nabla c|_2| c|_\infty| u|_2\leq C| u|^2_2+C|\nabla c|^2_2|c|^2_\infty,\\
I_3=&\int_{\mathbb{R}^3} \psi \cdot Q(c,v)\cdot u  \text{d}x\\
\leq&  C(1+|\overline{E}(c)|_{\infty})|\psi|_6|\nabla v|_3|u|_2\leq C|u|^2_2+M(b_0)|\psi|^2_6|\nabla v|^2_3.
\end{split}
\end{equation}
Then we have
\begin{equation}\label{zhu3}
\begin{split}
&\frac{1}{2} \frac{d}{dt}|u|^2_2+\alpha|\nabla u|^2_2\leq C(|u|^2_2+|v|^2_3|\nabla v|^2_2+|\nabla c|^2_2|c|^2_\infty)+M(b_0)|\psi|^2_6|\nabla v|^2_3.
\end{split}
\end{equation}
Integrating (\ref{zhu3}) over $(0,t)$, for $0 \leq t \leq T_1$, we have
\begin{equation*}
\begin{split}
|u(t)|^2_2+\int_0^t\alpha|\nabla u(s)|^2_2\text{d}s\leq  C\int_0^t  |u(s)|^2_2 \text{d}s+C|u_0|^2_2+M(b_0)b^{4}_3t.
\end{split}
\end{equation*}
According to Gronwall's inequality, we have
\begin{equation}\label{zhu5}
\begin{split}
|u(t)|^2_2+\int_0^t\alpha|\nabla u(s)|^2_2\text{d}s\leq  C\big(|u_0|^2_2+M(b_0)b^{4}_3t\big)\exp(Ct)\leq Cb^2_0
\end{split}
\end{equation}
for $0 \leq t \leq T_2=\min(T^*,(1+M(b_0)b^{4}_3)^{-1})$.
\end{proof}

Next, in order to obtain the higher order regularity estimate for the velocity $u$,  we need to  introduce the effective viscous flux $F$ and vorticity $\omega$ to deal with the $c$-dependent  Lam$\acute{ \text{e} }$ operator  (see (\ref{eq:5.231})), which can be given as
\begin{equation}\label{liu09}
\begin{split}
F=(2\alpha+ \overline{E}(c))\text{div}u-(\theta  c^2-\theta (c^{\infty})^2),\quad \omega=\nabla \times u,
\end{split}
\end{equation}
then in the sence of distribution,  the momentum equations $(\ref{li4})_2$ can be written as 
\begin{equation}\label{liu010}
\begin{cases}
\triangle F=\text{div}(u_t+v\cdot \nabla v-\psi \cdot Q(c,v)), \\[8pt]
 \triangle \omega=\nabla \times (u_t+v\cdot \nabla v-\psi \cdot Q(c,v)).
\end{cases}
\end{equation}
So  we immediately have
\begin{equation}\label{liu099}
\begin{split}
-\triangle u=\nabla \times \omega- \nabla \text{div}u=\nabla \times \omega- \nabla \Big(\frac{F+\theta c^2-\theta (c^{\infty})^2}{2\alpha+ \overline{E}(c)}\Big).
\end{split}
\end{equation}

 \begin{lemma}[\textbf{Higher order estimates of the velocity $u$}]\label{ssk4}
\begin{equation*}
\begin{split}
|u(t)|^2_{D^1}+\int_{0}^{t}\Big(|u_t(s)|^2_{2}+|u(s)|^2_{D^2}\Big)\text{d}s\leq& Cb^{2}_0,\\
| u(t)|^2_{D^2}+|u_t(t)|^2_{2}+\int_{0}^{t}\Big(|u(s)|^2_{D^{3}}+|u_t(s)|^2_{D^1}\Big)\text{d}s\leq& M(b_0)b^{3}_2b_3,
\end{split}
\end{equation*}
for $0 \leq t \leq T_3=\min(T^*,(1+M(b_0)b^{8}_3)^{-1})$.
 \end{lemma}

\begin{proof}
\underline{Step 1}. 
Via the standard elliptic estimate shown in Lemma \ref{zhenok} and (\ref{liu099}), we immediately obtain 
\begin{equation*}
\begin{split}
|u|_{ D^2}\leq& C (|\nabla \times \omega|_{ 2}+|\nabla F|_{ 2}+|\nabla c^2|_{ 2}+|\nabla \overline{E}(c)|_{ 6}|\text{div} u|_{ 3})\\
\leq &C (|\nabla \times \omega|_{ 2}+|\nabla F|_{ 2}+|\nabla c^2|_{ 2}+|\nabla u|_2|\nabla \overline{E}(c)|^2_{ 6})+\frac{1}{2}|u|_{ D^2},\\
\end{split}
\end{equation*}
where we have used the fact that 
\begin{equation}\label{heavn1}
\begin{split}
\text{div}u=\frac{F+\theta c^2-\theta (c^{\infty})^2}{2\alpha+ \overline{E}(c)},\ \text{and} \ |\text{div} u|_3\leq C|\nabla u|^{\frac{1}{2}}_2|\nabla u|^{\frac{1}{2}}_6.
\end{split}
\end{equation}
Then via Young's inequality,  we have
\begin{equation}\label{liu11}
\begin{split}
|u|_{ D^2}\leq C (M(b_0)|\nabla u|_2+|\nabla  \omega|_{ 2}+|\nabla F|_{ 2}+b^2_0).
\end{split}
\end{equation}
Again from Lemma \ref{zhenok}, we also have
\begin{equation}\label{liu12}
\begin{split}
|\nabla  \omega|_{ 2}+|\nabla F|_{ 2}\leq &C(|u_t|_{ 2}+|v|_{ 6}|\nabla v|_{ 3}+|\psi|_{ 6}|Q(c, v)|_{ 3})\leq C(M(b_0)b^{\frac{3}{2}}_2b^{\frac{1}{2}}_3+|u_t|_{ 2}).
\end{split}
\end{equation}
Then combining (\ref{liu11})-(\ref{liu12}), we deduce that 
\begin{equation}\label{liu13}
\begin{split}
|u|_{ D^2}\leq C (M(b_0)|\nabla u|_2+|u_t|_{ 2}+M(b_0)b^{\frac{3}{2}}_2b^{\frac{1}{2}}_3).
\end{split}
\end{equation}
\underline{Step 2} (Estimate for $| \nabla u|_2$). Multiplying $(\ref{li4})_2$ by $u_t$ and integrating over $\mathbb{R}^3$, we have
\begin{equation}\label{azhu8}
\begin{split}
&\frac{1}{2} \frac{d}{dt}\int_{\mathbb{R}^3}\big(\alpha |\nabla u|^2+(\alpha+ \overline{E}(c))|\text{div}u|^2\big)\text{d}x+|u_t|^2_2\\
=&\int_{\mathbb{R}^3} \Big(\frac{1}{2} \overline{E}(c)_t(\text{div}u)^2-\big((v\cdot \nabla v)+\theta( \nabla c^2) -(\psi \cdot Q(c,v))\big)\cdot u_t \Big) \text{d}x\equiv: \sum_{i=4}^7 I_i.
\end{split}
\end{equation}
According to H\"older's inequality, Lemma \ref{lem2as}, Young's inequality and (\ref{liu13}), 
\begin{equation}\label{azhu9}
\begin{split}
I_4=&\int_{\mathbb{R}^3} \frac{1}{2} \overline{E}(c)_t(\text{div}u)^2\text{d}x\leq C|\overline{E}(c)_t|_{3}|\nabla u|_{2}|\nabla u|_6\\
\leq& \epsilon |u|^2_{D^2}+C(\epsilon)|\overline{E}(c)_t|^2_{3}| u|^2_{D^1},\\
I_5=&-\int_{\mathbb{R}^3} (v\cdot \nabla v) \cdot u_t \text{d}x\leq C|v|_{\infty}|\nabla v|_{2}|u_t|_2\\
\leq& C\|\nabla v\|^2_{1}|\nabla v|^2_{2}+\frac{1}{10}|u_t|^2_2,\\
I_6=&-\int_{\mathbb{R}^3} 2\theta(c \nabla c) \cdot u_t \text{d}x
\leq C|\nabla c|_2| c|_\infty|u_t|_2\\
\leq & \frac{1}{10} |u_t|^2_2+C|\nabla c|^2_2| c|^2_\infty,\\
I_7=&\int_{\mathbb{R}^3} \psi \cdot Q(c,v)\cdot u_t\text{d}x\leq C|u_t|_2|\psi|_6|Q(c,v)|_3\\
\leq & \frac{1}{10}|u_t|^2_2+C|\psi|^2_6|Q(c,v)|^2_3,
\end{split}
\end{equation}
where $\epsilon>0$ is a sufficiently small constant.

Combining  (\ref{liu13}) and (\ref{azhu8})-(\ref{azhu9}),  via letting $\epsilon$ sufficiently small,  we have
\begin{equation}\label{zhu8qq}
\begin{split}
& \frac{d}{dt}|\nabla u|^2_2+| u_t|^2_2
\leq M(b_0)b^4_3|\nabla u|^2_2+M(b_0)b^4_3.
\end{split}
\end{equation}
From Gronwall's inequality, we have
\begin{equation}\label{zhu8qqfg}
\begin{split}
|\nabla u(t)|^2_2+\int_0^t |u_t|^2_2\text{d}s \leq C(|\nabla u_0|^2_2+M(b_0)b^4_3t)\exp(M(b_0)b^4_3t)\leq Cb^2_0,\end{split}
\end{equation}
for $0 \leq t \leq T'=\min(T^*,(1+M(b_0)b^{4}_3)^{-1})$, which, along with (\ref{liu13}), implies that 
$$
\int_0^t |u|^2_{ D^2}\leq C \int_0^t \Big(M(b_0)|\nabla u|_2+|u_t|_{ 2}+M(b_0)b^{\frac{3}{2}}_2b^{\frac{1}{2}}_3\Big)^2 \text{d}s\leq Cb^2_0,\quad \text{for}\quad 0 \leq t \leq T'.
$$

\underline{Step 3} (Estimate for $| \nabla^2 u|_2$).
 We consider the estimate for $| u_t|_2$. First we differential $(\ref{li4})_2$ with respect to $t$:
\begin{equation}\label{zhu7}
\begin{split}
u_{tt}+(Lu)_t=-(v\cdot\nabla v)_t -2\theta(c\nabla c)_t+(\psi\cdot Q(c,v))_t.
\end{split}
\end{equation}
Then multiplying (\ref{zhu7}) by $u_t$ and integrating over $\mathbb{R}^3$, we have
\begin{equation}\label{zhu8}
\begin{split}
&\frac{1}{2} \frac{d}{dt}|u_t|^2_2+\alpha|\nabla u_t|^2_2+\int_{\mathbb{R}^3}(\alpha+ \overline{E}(c))|\text{div} u_t|^2\text{d}x\\
=&\int_{\mathbb{R}^3} \big(- \overline{E}(c)_t\text{div}u\text{div}u_t-\big((v\cdot \nabla v)_t +\theta( \nabla c^2)_t -(\psi \cdot Q(c,v))_t\big)\cdot u_t \big) \text{d}x\equiv: \sum_{i=8}^{11} I_i.
\end{split}
\end{equation}
According to H\"older's inequality, Lemma \ref{lem2as} and Young's inequality, 
\begin{equation}\label{zhu9}
\begin{split}
I_8=&-\int_{\mathbb{R}^3}  \overline{E}(c)_t\text{div}u\text{div}u_t \text{d}x\leq C|\overline{E}(c)_t|_{3}|\nabla u_t|_{2}|\nabla u|_6\\
\leq& \frac{\alpha}{10}|\nabla u_t|^2_{2}+C|\overline{E}(c)_t|^2_{3}| u|^2_{D^2},\\
I_9=&-\int_{\mathbb{R}^3} (v\cdot \nabla v)_t \cdot u_t \text{d}x\leq C\big(|v|_{\infty}|\nabla v_t|_{2}|u_t|_2+|v_t|_{6}|\nabla v|_{3}|u_t|_2\big)\\
\leq& \frac{1}{b^2_3}|\nabla v_t|^2_{2}+Cb^2_3\|\nabla v\|^2_{1}|u_t|^2_2,\\
I_{10}=&-\int_{\mathbb{R}^3} \theta( \nabla c^2)_t \cdot u_t \text{d}x=\theta\int_{\mathbb{R}^3}   (c^2)_t \text{div}u_t \text{d}x\\
\leq& C|c_t|_2| c|_\infty|\nabla u_t|_2\leq \frac{\alpha}{10}|\nabla u_t|^2_2+C|c_t|^2_2|c|^2_\infty.
\end{split}
\end{equation}
For  the last term on the right side of (\ref{zhu8}),  we have
\begin{equation}\label{liu05}\begin{split}
I_{11}=&\int_{\mathbb{R}^3} \psi \cdot Q(c,v)_t\cdot u_t\text{d}x+\int_{\mathbb{R}^3} \psi_t \cdot Q(c,v)\cdot u_t\text{d}x
=I_{11A}+I_{11B}.
\end{split}
\end{equation}
We firstly consider the term:
\begin{equation}\label{liu06}\begin{split}
I_{11A}
\leq& C(1+|\overline{E}(c)|_{\infty} )|\psi|_6 |\nabla v_t|_2|u_t|_{3}+C|\overline{E} (c)_t|_{2} |\psi|_6|\nabla v|_{\infty}|u_t|_{3}\\
\leq&\frac{1}{b^2_3}|\nabla v_t|^2_{2}+\frac{\alpha}{10}|\nabla u_t|^2_2+M(b_0)b^8_3|u_t|^2_2+M(b_0)b^8_3,
\end{split}
\end{equation}
where we have used the fact that 
\begin{equation}\label{liu07}
|u_t|_{3}\leq C u_t|^{\frac{1}{2}}_{2}|\nabla u_t|^{\frac{1}{2}}_{2},\ |Q(c,v)_t|_{2}\leq C (1+|\overline{E}(c)|_{\infty})|\nabla v_t|_2+C|\overline{E}(c)_t|_{2}|\nabla v|_\infty.
\end{equation}
And for the second term:
\begin{equation}\label{liu07}\begin{split}
I_{11B}=&-\int_{\mathbb{R}^3}\Big(\nabla (v \cdot \psi) \cdot Q(c,v)\cdot u_t+\nabla \text{div}v\cdot Q(c,v)\cdot u_t\Big)\text{d}x\\
\leq& C\int_{\mathbb{R}^3}\big(|v| |\psi| |\nabla Q(c,v)| |u_t|+ |v| |\psi| |Q(c,v)| |\nabla u_t|+ |\nabla \text{div}v| |Q(c,v)| |u_t| \big)\text{d}x  \\
\leq &C | v|_{\infty} |\psi|_6\big( |\nabla Q(c,v)|_2|u_t|_{3}+ | Q(c,v)|_{3}|\nabla u_t|_{2}\big)+C |\nabla^2 v|_6 | Q(c,v)|_2 |u_t|_{3}\\
\leq&\frac{1}{b^2_3}| v|^2_{D^{2,6}}+\frac{\alpha}{10}|\nabla u_t|^2_2+M(b_0)b^8_3|u_t|^2_2+M(b_0)b^8_3,
\end{split}
\end{equation}
where we have used the fact that 
\begin{equation}\label{liu08}
\begin{cases}
|\nabla Q(c,v)|_{2}\leq C(1+|\overline{E}(c) |_{\infty})|v|_{D^2}+C|\nabla \overline{E}(c)|_{3}|\nabla v|_{6},\\[8pt]
|Q(c,v)|_{3}\leq C|Q(c,v) |^{\frac{1}{2}}_{2}|\nabla Q(c,v)|^{\frac{1}{2}}_{2},\ | Q(c,v)|_{2} \leq C(1+|\overline{E}(c) |_{\infty})|v|_{D^1}.
\end{cases}
\end{equation}
Combining (\ref{liu13}),  (\ref{zhu8qqfg}) and (\ref{zhu8})-(\ref{liu07}),  we have
\begin{equation}\label{zhu8qq}
\begin{split}
& \frac{1}{2}\frac{d}{dt}|u_t|^2_2+\alpha|\nabla u_t|^2_2+\int_{\mathbb{R}^3}(\alpha+\overline{ E}(c))|\text{div} u_t|^2 \text{d}x\\
\leq &M(b_0)b^8_3|u_t|^2_2+M(b_0)b^{8}_3+\frac{C}{b^2_3}\big(|\nabla v_t|^2_{2}+|v|^2_{D^{2,6}}\big).
\end{split}
\end{equation}

Integrating (\ref{zhu8qq}) over $(\tau,t)$ $(\tau \in( 0,t))$ for $ 0< t \leq T_3$, we have
\begin{equation}\label{zhu13}
\begin{split}
&|u_t(t)|^2_2+\int_\tau^t \alpha|\nabla u_t(s)|^2_2\text{d}s\\
\leq & |u_t(\tau)|^2_2+M(b_0)b^{8}_3t+\int_{\tau}^tM(b_0)b^8_3 |u_t|^2_2\text{d}s+C.
\end{split}
\end{equation}

According to the momentum equations $ (\ref{li4})_2$, we have
\begin{equation}\label{zhu15}
\begin{split}
|u_t(\tau)|_2\leq C\big( |v|_\infty |\nabla v|_2+|c|_\infty|\nabla c|_2+|Lu|_{2}+|\psi|_6|Q(c,v)|_2\big)(\tau).
\end{split}
\end{equation}
Then via the assumptions (\ref{vg})-(\ref{reggh}), we easily have
\begin{equation}\label{zhu15vb}
\begin{split}
\lim \sup_{\tau\rightarrow 0}|u_t(\tau)|_2\leq& C\big( |v_0|_\infty |\nabla v_0|_2+|c_0|_\infty|\nabla c_0|_2+|Lu_0|_{2}+|\psi_0|_6|Q(c_0, v_0)|_2\big)\\
\leq& M(b_0)b^2_0.
\end{split}
\end{equation}

So letting $\tau\rightarrow 0$ in (\ref{zhu13}), via Gronwall's inequality, we have
\begin{equation}\label{zhu14}
\begin{split}
&|u_t(t)|^2_2+\int_0^t \alpha|\nabla u_t(s)|^2_2\text{d}s
\leq \big(M(b_0)b^{8}_3t+M(b_0)b^2_0\big)\exp (M(b_0)b^8_3t)\leq M(b_0)b^2_0,
\end{split}
\end{equation}
for $0 \leq t \leq T_3=\min(T^*,(1+M(b_0)b^{8}_3)^{-1})$.

\underline{Step 4}. Finally, we consider the estimates of the higher order terms.
From  estimate  (\ref{liu13}), Lemmas \ref{lem2as} and  \ref{zhenok},  relation (\ref{heavn1}) and inequality (\ref{ine}),
 we easily have, for $0 \leq t \leq T_3$,
\begin{equation*}
\begin{split}
|u(t)|_{D^2}\leq C (M(b_0)|\nabla u|_2+|u_t|_{ 2}+M(b_0)b^{\frac{3}{2}}_2b^{\frac{1}{2}}_3)\leq M(b_0)b^{\frac{3}{2}}_2b^{\frac{1}{2}}_3,
\end{split}
\end{equation*}
and
\begin{equation}\label{ji1}
\begin{split}
|u|_{D^{3}}\leq& C \big(|\nabla \omega|_{D^1}+|\nabla F|_{D^1} +|\nabla c^2|_{D^1}+|\nabla \overline{E}(c)|_{3}|\nabla F|_6\big)\\
&+C\big(|\nabla \overline{E}(c)|_{3}(|\nabla c^2|_6+|\nabla u|_6|\nabla \overline{E}(c)|_2)+|\nabla u|_\infty| \overline{E}(c)|_{D^2}\big)\\
\leq & M(b_0) \big(|\nabla \omega|_{D^1}+|\nabla F|_{D^1} + b^{\frac{3}{2}}_2b^{\frac{1}{2}}_3)+\frac{1}{2}|u|_{D^3}.
\end{split}
\end{equation}
Again from Lemma \ref{zhenok} and (\ref{liu010}), we also have
\begin{equation*}
\begin{split}
|\nabla \omega|_{D^1}+|\nabla F|_{D^1} \leq& C(|u_t|_{ D^1}+|v \cdot \nabla v|_{ D^1}+|\psi \cdot Q(c,v)|_{D^1})\\
\leq &C(|u_t|_{ D^1}+M(b_0)b^3_3),
\end{split}
\end{equation*}
which, together with (\ref{zhu14})-(\ref{ji1}), immediately implies the desired estimate for $|u|_{D^3}$.
\end{proof}

Then combining the estimates obtained in Lemmas \ref{2}-\ref{ssk4}, we have
\begin{equation}\label{jkkas}
\begin{split}
|c(t)|^2_\infty+\|c(t)-c^{\infty}\|^2_2+\|c_t(t)\|^2_1\leq M(b_0)b^4_3,\\
|\overline{E}(c)(t)|^2_\infty+\|\overline{E}(c)(t)-\overline{E}(c^\infty)\|^2_{2}+\|\overline{E}(c)_t(t)\|^2_{1}\leq M(b_0)b^4_3,\\
|\psi(t)|^2_{D^1}+|\psi(t)_t|^2_{ 2}\leq M(b_0)b^4_3,\\
\|u(t)\|^2_{1}+\int_{0}^{t}\Big(|u_t(s)|^2_{2}+\|\nabla u(s)\|^2_{1}\Big)\text{d}s\leq M(b_0)b^{2}_0,\\
| u(t)|^2_{D^2}+|u_t(t)|^2_{2}+\int_{0}^{t}\Big(|u(s)|^2_{D^{3}}+|u_t(s)|^2_{D^1}\Big)\text{d}s\leq M(b_0)b^{3}_2b_3
\end{split}
\end{equation}
for $0 \leq t \leq T_3$. Therefore, if we define the constants $b_i$ ($i=1,2,3$) and $T^*$ by
\begin{equation}\label{dingyi}
\begin{split}
&b_1=b_2=M(b_0)b_0,\ b_3= M(b_0)b^3_2=M^4(b_0)b^3_0, \\
& \text{and} \ T^*=\min (T, (1+M(b_0)b_3)^{-8}),
\end{split}
\end{equation}
then we deduce that 
\begin{equation}\label{jkk}
\begin{split}
\sup_{0\leq t \leq T^*}| u(t)|^2_{2}+\int_{0}^{T^*} |\nabla u(t)|^2_{2}\text{d}t \leq b^2_1&,\\
\sup_{0\leq t \leq T^*}| u(t)|^2_{D^1}+\int_{0}^{T^*} \Big( | u(t)|^2_{D^2}+|u_t(t)|^2_{2}\Big)\text{d}t \leq b^2_2&,\\
\sup_{0\leq t \leq T^*}(|u(t)|^2_{D^2}+|u_t(t)|^2_{2})+\int_{0}^{T^*} \Big( |u(t)|^2_{D^3}+|u_t(t)|^2_{D^1}\Big)\text{d}t \leq b^2_3&,\\
\sup_{0\leq t \leq T^*}\big(|c(t)|^2_\infty+\|c(t)-c^{\infty}\|^2_2+\|c_t(t)\|^2_1\big)\leq M(b_0)b^4_3&,\\
\sup_{0\leq t \leq T^*}\big(|\overline{E}(c)(t)|^2_\infty+\|\overline{E}(c)(t)-\overline{E}(c^\infty)\|^2_{2}+\|\overline{E}(c)_t(t)\|^2_{1}\big)\leq M(b_0)b^4_3&,\\
\sup_{0\leq t \leq T^*}\big(|\psi(t)|^2_{D^1}+|\psi(t)_t|^2_{ 2}\big)\leq M(b_0)b^4_3&.
\end{split}
\end{equation}

\subsection{Unique solvability of the linearization  with vacuum}
Based on the  a prior estimate (\ref{jkk}), we have the following existence result under the assumption that $c_0\geq 0$.
\begin{lemma}\label{lem1q}
 Assume that  the initial data (\ref{qwe1}) satisfy (\ref{th78rr}) and $c_0\geq 0$.
Then there exists a unique strong solution $(c, \psi ,u)$ to (\ref{li4})-(\ref{vg}) such that
\begin{equation}\label{regghq}\begin{split}
&c\geq 0, \ c \in C([0,T^*]; H^2), \ c _t \in C([0,T^*]; H^1),\ \psi \in C([0,T^*] ; D^1), \\
&\psi_t \in C([0,T^*]; L^2),\  u\in C([0,T^*]; H^2)\cap L^2([0,T^*] ; H^3), \\
 & u_t \in C([0,T^*]; L^2)\cap L^2([0,T^*] ; D^1).
\end{split}
\end{equation}
And we also have 
$\partial_i \psi^{(j)}=\partial_j \psi^{(i)}$ in the distribution sense for 
 $(i,j=1,2,3)$.
Moreover,  $(c, \psi ,u)$ also satisfies the local a prior estimates (\ref{jkk}).
\end{lemma}
\begin{proof}
\underline{Step 1}. Existence.
We firstly define 
$$c_{\delta0}=c_0+\delta,\quad \text{and} \quad \psi_{\delta0}=2\theta\nabla c_0/(c_0+\delta)$$
 for each $\delta\in (0,1)$. 
Then according to the assumption (\ref{houmian}), for all sufficiently small $\delta> 0$,
\begin{equation*}\begin{split}
1+|c_{\delta0}|_\infty+\|c_{\delta0}-\delta\|_{2}+|\psi_{\delta0}|_{D^1}+\|u_0\|_{2}\leq Cb^2_0= \overline{b}_0.
\end{split}
\end{equation*}
Therefore, corresponding to $(\rho_{\delta0},u_0,\psi_{\delta0})$ with small $\delta> 0$, there exists a unique strong solution 
$(c^\delta,u^\delta,\psi^\delta)$  to the linearized problem (\ref{li4})-(\ref{vg}) satisfying the local estimate (\ref{jkk}) obtained in the above section. 

By virtue of this uniform estimates (\ref{jkk}), we know that  there exists a subsequence of solutions 
\begin{equation}\label{weak}\begin{split}
(c^\delta,u^\delta,\psi^\delta)\quad  \text{converges \  to \  a \  limit} \quad (c,u,\psi) \quad \text{ in \  weak \ or \ weak*\ sense}.
\end{split}
\end{equation}
 And for any $R> 0$, due to the compact property in Lemma \ref{aubin} (see \cite{jm}), there exists a subsequence of solutions $(c^\delta,u^\delta,\psi^\delta)$ satisfy:
\begin{equation}\label{ert}\begin{split}
(c^\delta,u^\delta)\rightarrow (c,u) \ \text{in } \ C([0,T^*];H^1(B_R)),\quad \psi^\delta \rightarrow \psi \ \text{in } \ C([0,T^*];L^2(B_R)),
\end{split}
\end{equation}
where $B_R$ is a ball centered at origin with radius $R$.
Combining the lower semi-continuity of norms, the  weak or weak* convergence of $(c^\delta,u^\delta,\psi^\delta)$ and (\ref{ert}), we know that $(c,u,\psi) $ also satisfies the local estimates (\ref{jkk}). 

Then via the local estimates (\ref{jkk}),   the weak or weak* convergence in (\ref{weak}) and strong convergence in (\ref{ert}), in order to make sure that  $(c,u,\psi) $ is a weak solution in  the sense of distribution to the linearized problem (\ref{li4})-(\ref{vg})
satisfying the  regularity
\begin{equation}\label{wode}\begin{split}
&c\geq 0, \ c \in L^\infty([0,T^*]; H^2), \ c _t \in L^\infty([0,T^*]; H^1), \\
& \psi \in L^\infty([0,T^*] ; D^1),\  \psi_t \in L^\infty([0,T^*]; L^2), \ u\in L^\infty([0,T^*]; H^2)\cap L^2([0,T^*] ; H^3), \\
 &u_t \in L^\infty([0,T^*]; L^2)\cap L^2([0,T^*] ; D^1),
\end{split}
\end{equation} we only need to make sure that 
\begin{equation}\label{initial}\begin{split}
&\lim_{\delta\rightarrow 0} \int_{\mathbb{R}^3} (c^\delta_0-c^0) \phi(0,x) \text{d}x=0,\\
&\lim_{\delta\rightarrow 0} \int_{\mathbb{R}^3} (\psi^\delta_0-\psi^0) \xi(0,x) \text{d}x=0
\end{split}
\end{equation}
for any $\phi(t,x)\in C^\infty_c([0,T^*)\times \mathbb{R}^3)$ and $\xi(t,x) \in C^\infty_c([0,T^*)\times \mathbb{R}^3)^3$. The proof for $(\ref{initial})_1$ is easy, so we only need to consider $(\ref{initial})_2$. When
$$\supp_{x}\xi (0,x) \cap \{x\in \mathbb{R}^3|c_0(x)=0\}=\emptyset,$$
then due to $c_0\in H^2( \mathbb{R}^3) \subset C( \mathbb{R}^3)$, there must exists a positive constant $\delta_0$ such that 
\begin{equation}\label{zheng1}
c_0(x)>\delta_0 \quad \text{for} \quad x\in \supp_{x}\xi (0,x),
\end{equation}
which immedaitely implies that 
\begin{equation}\label{chu1}\begin{split}
&\lim_{\delta\rightarrow 0} \int_{\mathbb{R}^3} (\psi^\delta_0-\psi^0) \xi(0,x) \text{d}x
=\lim_{\delta\rightarrow 0} \int_{\supp_{x}\xi (0,x)} -\frac{\delta}{c_0+\delta}\psi^0 \xi(0,x) \text{d}x\\
\leq & \lim_{\delta\rightarrow 0}\frac{\delta}{\delta_0+\delta} |\xi (0,x)|_2 |\psi_0|_6|\supp_{x}\xi (0,x)|^{\frac{1}{3}}\rightarrow 0,
\end{split}
\end{equation}
where $|\supp_{x}\xi (0,x)|$ means the $3$D Lebesgue measure of $\supp_{x}\xi (0,x)$.

And when
$$\supp_{x}\xi (0,x) \cap \{x\in \mathbb{R}^3|c_0(x)=0\}\neq \emptyset,$$
due to $\psi_0=\nabla c_0 / c_0 \in D^1(\mathbb{R}^3)$, we must have 
$$
|\{x\in \mathbb{R}^3|c_0(x)=0\}|=0.
$$
Then for every $n\geq 1$, we have  
\begin{equation}\label{chu2}\begin{split}
&I= \int_{\mathbb{R}^3} (\psi^\delta_0-\psi^0) \xi(0,x) \text{d}x
=\int_{\supp_{x}\xi (0,x)} -\frac{\delta}{c_0+\delta}\psi^0 \xi(0,x) \text{d}x\\
=& \int_{\supp_{x}\xi (0,x)\cap \{x\in \mathbb{R}^3|c_0(x)\geq \frac{1}{n}\} } -\frac{\delta}{c_0+\delta}\psi^0 \xi(0,x) \text{d}x\\
&+ \int_{\supp_{x}\xi (0,x)\cap \{x\in \mathbb{R}^3|c_0(x)<\frac{1}{n}\} } -\frac{\delta}{c_0+\delta}\psi^0 \xi(0,x) \text{d}x=I_1+I_2.
\end{split}
\end{equation}
So it is easy to see that 
\begin{equation}\label{chu3}\begin{split}
\lim_{\delta \rightarrow 0}I=&\lim_{n \rightarrow +\infty}\lim_{\delta \rightarrow 0}I=\lim_{n \rightarrow +\infty}\lim_{\delta \rightarrow 0}I_2\\
\leq& C |\xi (0,x)|_2 |\psi_0|_6 \lim_{n \rightarrow +\infty}|\supp_{x}\xi (0,x)\cap \{x\in \mathbb{R}^3|c_0(x)<1/n\} |^{\frac{1}{3}}=0,
\end{split}
\end{equation}
which, together with (\ref{chu1}), implies that (\ref{initial}) holds. 

Moreover, from the conclusions obtained in this step,  we also know that  even vacuum appears, $\psi$ satisfies $\partial_i \psi^{(j)}=\partial_j \psi^{(i)} \ (i,j=1,2,3)$  and  the following positive and symmetric hyperbolic system in the distribution sense:
\begin{equation}\label{zhenzheng}
\psi_t+\sum_{l=1}^3 A_l \partial_l\psi+B\psi+\nabla \text{div} v=0, \ \psi_0 \in  D^1.
\end{equation}

\underline{Step 2}. The uniqueness and time continuity for $(c,\psi,u)$ can be obtained via the same arguments used  in Lemma \ref{lem1}.

\end{proof}

\subsection{Proof of Theorem \ref{th2}} Our proof  is based on the classical iteration scheme and the existence results obtained  in  Section 3.3. Let us denote as in Section 3.2 that
\begin{equation*}\begin{split}
&2+|c_0|_{\infty}+\|(c_0,u_0)\|_{2}+|\psi_0|_{D^1} \leq b_0.
\end{split}
\end{equation*}
Next, let $u^0\in C([0,T^*];H^2)\cap  L^2([0,T^*];H^3) $ be the solution to the linear parabolic problem
$$
h_t-\triangle h=0 \quad \text{in} \quad (0,+\infty)\times \mathbb{R}^3 \quad \text{and} \quad h(0)=u_0 \quad \text{in} \quad \mathbb{R}^3.
$$
Then taking a small time $T^\epsilon\in (0,T^*]$, we have
\begin{equation}\label{jizhu}
\begin{split}
\sup_{0\leq t \leq T^\epsilon}| u^0(t)|^2_{2}+\int_{0}^{T^\epsilon} |\nabla u^0(t)|^2_{2}\text{d}t \leq b^2_1&,\\
\sup_{0\leq t \leq T^\epsilon}| u^0(t)|^2_{D^1}+\int_{0}^{T^\epsilon} \Big( |u^0(t)|^2_{D^2}+|u^0_t(t)|^2_{2}\Big)\text{d}t \leq b^2_2&,\\
\sup_{0\leq t \leq T^\epsilon}(|u^0(t)|^2_{D^2}+|u^0_t(t)|^2_{2})+\int_{0}^{T^\epsilon} \Big( |u^0(t)|^2_{D^3}+|u^0_t(t)|^2_{D^1}\Big)\text{d}t \leq b^2_3&.
\end{split}
\end{equation}

\begin{proof}
\underline{Step 1}. Existence.
Let $v=u^0$,  we can get $(c^1, \psi^1, u^1)$ as a strong solution to problem (\ref{li4})-(\ref{vg}). Then we construct approximate solutions $(c^{k+1}, \psi^{k+1}, u^{k+1})$ inductively, as follows: assuming that $(c^{k},\psi^k, u^{k})$ was defined for $k\geq 1$, let $(c^{k+1}, \psi^{k+1}, u^{k+1})$  be the unique solution to problem  (\ref{li4})-(\ref{vg}) with $v$ replaced by $ u^{k}$ as following:

\begin{equation}\label{li6}
\begin{cases}
c^{k+1}_t+u^{k}\cdot \nabla c^{k+1}+\frac{\gamma-1}{2}c^{k+1}\text{div} u^{k}=0,\\[8pt]
\psi^{k+1}_t+\sum_{l=1}^3 A_l(u^k) \partial_l\psi^{k+1}+B(u^k)\psi^{k+1}+\nabla \text{div}u^k=0,\\[8pt]
u^{k+1}_t+u^{k}\cdot\nabla u^{k} +2\theta c^{k+1}\nabla c^{k+1}=-L(c^{k+1})u^{k+1}+\psi^{k+1}\cdot Q(c^{k+1}, u^{k}),\\[8pt]
(c^{k+1}, \psi^{k+1}, u^{k+1})|_{t=0}=(c_0,\psi_0,u_0),\quad x\in \mathbb{R}^3,
 \end{cases}
\end{equation}
where the operator $L(f)g$ is defined as $
L(f)g=\text{div}(\alpha(\nabla g+(\nabla g)^\top)+ \overline{E}(f)\text{div}g\mathbb{I}_3)
$.
Via the estimates shown in Section $3.3$, we quickly deduce that the sequences of solutions $(c^{k}, \psi^{k}, u^{k})$ satisfy the uniform a prior estimate (\ref{jkk}).

The next task is to prove the strong convergence of the full sequence  $(c^k,\psi^k, u^k)$ of approximate solutions to a limit $(c,\psi,u)$ satisfying  (\ref{reg11}) in the sense of $H^1$.
Let
\begin{equation*}
\overline{c}^{k+1}=c^{k+1}-c^k,\ \overline{\psi}^{k+1}=\psi^{k+1}-\psi^k,\ \overline{u}^{k+1}=u^{k+1}-u^k,
\end{equation*}
 then from (\ref{li6}), we have
 \begin{equation}
\label{eq:1.2w}
\begin{cases}
\ \  \overline{c}^{k+1}_t+u^k\cdot \nabla\overline{c}^{k+1} +\overline{u}^k\cdot\nabla c ^{k}+\frac{\gamma-1}{2}(\overline{c}^{k+1} \text{div}u^k +c ^{k}\text{div}\overline{u}^k)=0,\\[8pt]
\ \ \overline{\psi}^{k+1}_t+\sum_{l=1}^3 A_l(u^k) \partial_l\overline{\psi}^{k+1}+B(u^k)\overline{\psi}^{k+1}+\nabla \text{div}\overline{u}^k=\Upsilon^k_1+\Upsilon^k_2,\\[8pt]
\ \ \overline{u}^{k+1}_t+ u^k\cdot\nabla \overline{u}^{k}+ \overline{u}^{k} \cdot \nabla u^{k-1}+\theta\nabla ((c^{k+1})^2-(c^k)^2)+L(c^{k+1})\overline{u}^{k+1} \\[8pt]
=\text{div}((\overline{E}(c^{k+1})-\overline{E}(c^{k}) ) \text{div}u^k \mathbb{I}_3)+\psi^{k+1}\cdot Q(c^{k+1},\overline{u}^k)\\[8pt]
\ \ +\overline{\psi}^{k+1}\cdot Q(c^k,u^{k-1})+\psi^{k+1}(\overline{E}(c^{k+1})-\overline{E}(c^{k}) ) \text{div}u^{k-1},
\end{cases}
\end{equation}
where $\Upsilon^k_1$  and $\Upsilon^k_2$ are defined via
\begin{equation*}
\Upsilon^k_1=-\sum_{l=1}^3(A_l(u^k) \partial_l\psi^{k}-A_l(u^{k-1}) \partial_l\psi^{k}),\quad \Upsilon^k_2=-(B(u^k) \psi^{k}-B(u^{k-1}) \psi^{k}).
\end{equation*}

Firstly multiplying $(\ref{eq:1.2w})_1$ by $2\overline{c}^{k+1}$ and integrating  over $\mathbb{R}^3$, we have
\begin{equation*}
\begin{split}
\frac{d}{dt}|\overline{c}^{k+1}|^2_2=& -2\int_{\mathbb{R}^3}\Big(u^k\cdot \nabla\overline{c}^{k+1} +\overline{u}^k\cdot\nabla c ^{k}+\frac{\gamma-1}{2}(\overline{c}^{k+1} \text{div}u^k +c ^{k}\text{div}\overline{u}^k)\Big)\overline{c}^{k+1} \text{d}x,\\
\leq& C|\nabla u^k|_\infty|\overline{c}^{k+1}|^2_2+C |\overline{c}^{k+1}|_2|\overline{u}^k|_6|\nabla c^k|_3+C|\overline{c}^{k+1}|_2|\nabla\overline{u}^k|_2| c^k|_\infty,
\end{split}
\end{equation*}
which means that ($0<\eta \leq \min\Big(\frac{1}{10},\frac{\alpha}{10}\Big)$ is a constant)
\begin{equation}\label{go64}\begin{cases}
\displaystyle
\frac{d}{dt}|\overline{c}^{k+1}(t)|^2_2\leq A^k_\eta(t)|\overline{c}^{k+1}(t)|^2_2+\eta |\nabla\overline{u}^k(t)|^2_2,\\[10pt]
\displaystyle
A^k_\eta(t)=C\Big(\|\nabla u^k\|_{2}+\frac{1}{\eta}\|c^{k}\|^2_{2} \Big),\ \text{and} \ \int_0^t A^k_\eta(s)\text{d}s\leq \widehat{C}+\widehat{C}_{\eta}t
\end{cases}
\end{equation}
for $t\in[0,T^\epsilon]$, where $\widehat{C}_{\eta}$ is a positive constant  depending on $\eta$ and  constant $\widehat{C}$.

Next, differentiating $(\ref{eq:1.2w})_1$ $\zeta$-times ($|\zeta|=1$) with respect to $x$, multiplying the resulting equation by $2 D^\zeta\overline{c}^{k+1}$ and  integrating  over $\mathbb{R}^3$,  we have
\begin{equation*}
\begin{split}
\frac{d}{dt}|D^\zeta\overline{c}^{k+1}|^2_{2}=&-2 \int_{\mathbb{R}^3}D^\zeta\big(u^k\cdot \nabla\overline{c}^{k+1} +\overline{u}^k\cdot\nabla c ^{k}+\frac{\gamma-1}{2}(\overline{c}^{k+1} \text{div}u^k +c ^{k}\text{div}\overline{u}^k)\big)D^\zeta\overline{c}^{k+1}\text{d}x\\
\leq& C|\nabla u^k|_\infty |\nabla \overline{c}^{k+1}|^2_{2}+C|\nabla c^k|_{\infty}| \nabla \overline{u}^k|_2|\nabla \overline{c}^{k+1}|_2\\
&+C|\nabla \overline{c}^{k+1}|_{3}| \nabla \overline{u}^k|_6|\nabla^2 c^{k}|_2+C|\nabla^2 u^k|_3 |\nabla \overline{c}^{k+1}|^2_{2}\\
&+C|\nabla \overline{u}^k|_6 |\nabla \overline{c}^{k+1}|_{2}|\nabla c^{k}|_3+C| c^k|_{\infty}  |\nabla \text{div} \overline{u}^k|_2   |\nabla \overline{c}^{k+1}|_{2},
\end{split}
\end{equation*}
which means that
\begin{equation}\label{go64qqqq}
\begin{cases}
\displaystyle
\frac{d}{dt}|\nabla \overline{c}^{k+1}(t)|^2_2\leq B^k_\eta(t)|\nabla\overline{c}^{k+1}(t)|^2_2+\eta|\nabla \text{div} \overline{u}^k(t)|^2_2+\eta |\nabla\overline{u}^k(t)|^2_2,\\[10pt]
\displaystyle
B^k_\eta(t)=C\Big(\|\nabla u^k\|_{2}+\frac{1}{\eta}\|c^{k}\|^2_{2} \Big),\ \text{and} \ \int_0^t B^k_\eta(s)\text{d}s\leq \widehat{C}+\widehat{C}_{\eta}t
\end{cases}
\end{equation}
for $t\in[0,T^\epsilon]$. Then combining (\ref{go64})-(\ref{go64qqqq}), we easily have
\begin{equation}\label{go64qqqqqq}
\begin{cases}
\displaystyle
\frac{d}{dt}\| \overline{c}^{k+1}(t)\|^2_1\leq \Phi^k_\eta(t)\| \overline{c}^{k+1}(t)\|^2_1+\eta \|\nabla\overline{u}^k(t)\|^2_1,\\[8pt]
\displaystyle
 \int_0^t \Phi^k_\eta(s)\text{d}s\leq \widehat{C}+\widehat{C}_{\eta}t\quad \text{for}\quad  t\in[0,T^\epsilon].
\end{cases}
\end{equation}

Secondly, multiplying  $(\ref{eq:1.2w})_2$ by $2\overline{\psi}^{k+1}$ and integrating  over $\mathbb{R}^3$, we have
\begin{equation}\label{go64aa}
\begin{split}
\frac{d}{dt}|\overline{\psi}^{k+1}|^2_2\leq& C\Big(\sum_{l=1}^{3}|\partial_{l}A_l(u^k)|_\infty+|B(u^k)|_\infty\Big)|\overline{\psi}^{k+1}|^2_2\\
&+C(|\Upsilon^k_1 |_2+|\Upsilon^k_2|_2+|\nabla^2 \overline{u}^k|_2)|\overline{\psi}^{k+1}|_2.
\end{split}
\end{equation}
From H\"older's inequality, it is easy to deduce that
\begin{equation}\label{go64aa1}
\begin{split}|\Upsilon^k_1 |_2\leq C|\nabla \psi^{k}|_2 |\overline{u}^k|_\infty, \quad  |\Upsilon^k_2 |_2\leq C| \psi^{k}|_6 | \nabla \overline{u}^k|_3.
\end{split}
\end{equation}
From (\ref{go64aa})-(\ref{go64aa1}), for $t\in[0,T^\epsilon]$, we have
\begin{equation}\label{go64aa2}\begin{cases}
\displaystyle
\frac{d}{dt}|\overline{\psi}^{k+1}(t)|^2_2\leq \Psi^k_\eta(t)|\overline{\psi}^{k+1}(t)|^2_2+\eta ||\nabla\overline{u}^k(t)||^2_1,\\[8pt]
\displaystyle
\Psi^k_\eta(t)=C\Big(\|\nabla u^k\|_{2} +\frac{1}{\eta}\ |\psi^{k} |^2_{D^1}+\frac{1}{\eta}\Big),\ \text{and} \ \int_0^t \Psi^k_\eta(s)\text{d}s\leq \widehat{C}+\widehat{C}_{\eta}t.
\end{cases}
\end{equation}

Thirdly, multiplying $(\ref{eq:1.2w})_3$ by $2\overline{u}^{k+1}$ and integrating  over $\mathbb{R}^3$, we have
\begin{equation*}\begin{split}
&\frac{d}{dt}|\overline{u}^{k+1}|^2_2+2\alpha|\nabla\overline{u}^{k+1} |^2_2+\int_{\mathbb{R}^3}(\alpha+\overline{E}(c^{k+1}))|\text{div}\overline{u}^{k+1} |^2 \text{d}x\\
=& -2\int_{\mathbb{R}^3}\Big(-\text{div}\big((\overline{E}(c^{k+1})-\overline{E}(c^{k}) ) \text{div}u^k\mathbb{I}_3\big) + u^k \cdot \nabla \overline{u}^{k}+\overline{u}^{k} \cdot \nabla u^{k-1}\Big) \cdot \overline{u}^{k+1}\text{d}x\\
&-2\int_{\mathbb{R}^3}\Big(\theta\nabla \big((c^{k+1})^2-(c^k)^2\big)-\psi^{k+1}\cdot Q(c^{k+1},\overline{u}^k)\Big) \cdot \overline{u}^{k+1}\text{d}x\\
&+2\int_{\mathbb{R}^3}\Big(\overline{\psi}^{k+1}\cdot Q(c^k,u^{k-1})+\psi^{k+1}(\overline{E}(c^{k+1})-\overline{E}(c^{k}) ) \text{div}u^{k-1} \Big) \cdot \overline{u}^{k+1}\text{d}x\\
\leq &C|\overline{c}^{k+1}|_6|\text{div}u^k|_3|\nabla \overline{u}^{k+1}|_2+C|u^k|_\infty |\nabla\overline{u}^{k}|_2 |\overline{u}^{k+1}|_2\\
&+C|\overline{u}^{k}|_6 |\overline{u}^{k+1}|_2 |\nabla u^{k-1}|_3+C \big(|c^{k+1}|_\infty+|c^{k}|_\infty\big)|\nabla \overline{u}^{k+1}|_2 |\overline{c}^{k+1}|_2\\
& +C(1+|\overline{E}(c)|_\infty)|\psi^{k+1}|_6 |\nabla\overline{u}^{k}|_2 |\overline{u}^{k+1}|_3 +C|\overline{\psi}^{k+1}|_2 |Q(c^k, u^{k-1})|_\infty |\overline{u}^{k+1}|_2\\
&+C|\psi^{k+1}|_6|\overline{E}(c^{k+1})-\overline{E}(c^{k}) |_2 |\text{div}u^{k-1}|_6 |\overline{u}^{k+1}|_6,
\end{split}
\end{equation*}
which implies that
\begin{equation}\label{gogo1}\begin{split}
&\frac{d}{dt}|\overline{u}^{k+1}|^2_2+\alpha |\nabla\overline{u}^{k+1} |^2_2\\
\leq& E^k_\eta(t)\|\overline{u}^{k+1}\|^2_1+E^k_2(t)\|\overline{c}^{k+1}\|^2_{1}
+E^k_3(t)|\overline{\psi}^{k+1}|^2_{2}+\eta|\nabla\overline{u}^{k}|^2_2,
\end{split}
\end{equation}
where
\begin{equation*}
\begin{cases}
\displaystyle
E^k_\eta(t)=C\Big(1+\frac{1}{\eta}|u^{k}|^2_{\infty}+\frac{1}{\eta}| \nabla u^{k-1}|^2_{3}+\frac{1}{\eta}|\psi^{k+1}|^2_{6}\Big), \\[12pt] 
\displaystyle
E^k_2(t)=C\Big(|c^{k+1}|_\infty+|c^{k}|_\infty+|\text{div}u^k|_3+|\psi^{k+1}|_6|\text{div}u^{k-1}|_6\Big)^2,\\[12pt]
\displaystyle
E^k_3(t)=C|\nabla u^{k-1}|^2_\infty,
\end{cases}
\end{equation*}
and we also have  $$\int_0^t \big(E^k_\eta(s)+E^k_2(s)+E^k_3(s)\big)\text{d}s\leq \widehat{C}+\widehat{C}_\eta t$$
 for $t\in[0,T^\epsilon]$.

Next, differentiating  $(\ref{eq:1.2w})_3$ $\zeta$-times ($|\zeta|=1$) with respect to $x$, multiplying the resulting equation by $D^\zeta\overline{u}^{k+1}$ and  integrating  over $\mathbb{R}^3$,  we have
\begin{equation*}
\begin{split}
&\frac{1}{2}\frac{d}{dt}|D^\zeta\overline{u}^{k+1}|^2_{2}+\alpha|\nabla D^\zeta\overline{u}^{k+1} |^2_2+\int_{\mathbb{R}^3}(\alpha+\overline{E}(c^{k+1}))|D^{\zeta}\text{div}\overline{u}^{k+1} |^2 \text{d}x\\
=& \int_{\mathbb{R}^3} \Big(\text{div}(D^\zeta\overline{E}(c^{k+1}) \text{div}\overline{u}^{k+1}\mathbb{I}_3)+D^{\zeta}\text{div}((\overline{E}(c^{k+1})-\overline{E}(c^{k}) ) \text{div}u^k \mathbb{I}_3)\Big)\cdot D^\zeta\overline{u}^{k+1} \text{d}x \\
& \int_{\mathbb{R}^3} D^\zeta\big(-u^k\cdot\nabla \overline{u}^{k}-\overline{u}^{k} \cdot \nabla u^{k-1} \big)\cdot D^\zeta\overline{u}^{k+1} \text{d}x \\
&+\int_{\mathbb{R}^3} D^\zeta\Big(-\theta \nabla ((c^{k+1})^2-(c^k)^2)+\psi^{k+1}\cdot Q(c^{k+1},\overline{u}^k)\Big)\cdot D^\zeta\overline{u}^{k+1} \text{d}x
\\
&+\int_{\mathbb{R}^3} D^\zeta\Big(\overline{\psi}^{k+1}\cdot Q(c^k,u^{k-1})+\psi^{k+1}(\overline{E}(c^{k+1})-\overline{E}(c^{k}) ) \text{div}u^{k-1}\Big)\cdot D^\zeta\overline{u}^{k+1} \text{d}x=\sum_{i=1}^7 J_i.
\end{split}
\end{equation*}
Then from integration by parts, Lemma \ref{lem2as} and H\"older's inequality,
\begin{equation}
\begin{split}\label{litong3}
J_1=&\int_{\mathbb{R}^3} \text{div}(D^\zeta\overline{E}(c^{k+1}) \text{div}\overline{u}^{k+1}\mathbb{I}_3)\cdot D^\zeta\overline{u}^{k+1} \text{d}x\\
\leq& C|\nabla \overline{u}^{k+1}|_3 |\nabla^2\overline{u}^{k+1}|_2 |D^\zeta\overline{E}(c^{k+1})|_6\leq C|\nabla \overline{u}^{k+1}|^{\frac{1}{2}}_2 |\nabla^2\overline{u}^{k+1}|^{\frac{3}{2}}_2 |D^\zeta c^{k+1}|_6,\\
J_2=&\int_{\mathbb{R}^3} D^{\zeta}\text{div}((\overline{E}(c^{k+1})-\overline{E}(c^{k}) ) \text{div}u^k \mathbb{I}_3)\cdot D^\zeta\overline{u}^{k+1} \text{d}x\\
\leq & C|\nabla \overline{c}^{k+1}|_2|\text{div}u^k|_\infty|\nabla^2 \overline{u}^{k+1}|_2+C|\overline{c}^{k+1}|_6|\nabla \text{div}u^k|_3|\nabla^2 \overline{u}^{k+1}|_2,\\
J_3=&\int_{\mathbb{R}^3} -D^\zeta( u^k\cdot\nabla \overline{u}^{k})\cdot D^\zeta\overline{u}^{k+1} \text{d}x\\
\leq &C|\nabla u^k|_6 |\nabla\overline{u}^{k}|_2 |\nabla \overline{u}^{k+1}|_3+C| u^k|_\infty |\overline{u}^{k}|_{D^2} |\nabla \overline{u}^{k+1}|_2,\\
J_4=& \int_{\mathbb{R}^3} -D^\zeta(\overline{u}^{k} \cdot \nabla u^{k-1})\cdot D^\zeta\overline{u}^{k+1} \text{d}x\\
\leq &C|\nabla\overline{u}^{k}|_2 |\nabla\overline{u}^{k+1}|_3 |\nabla u^{k-1}|_6 +C|\overline{u}^{k}|_6 |\nabla \overline{u}^{k+1}|_3 |\nabla^2 u^{k-1}|_2,\\
J_5=&\int_{\mathbb{R}^3} -\theta D^\zeta(\nabla ((c^{k+1})^2-(c^k)^2) )\cdot D^\zeta\overline{u}^{k+1} \text{d}x \\
 \leq &C |\nabla c^{k+1}+\nabla c^{k}|_3|\nabla^2\overline{u}^{k+1}|_2 |\overline{c}^{k+1}|_6+C |(c^{k+1}+c^{k})|_\infty|\nabla^2\overline{u}^{k+1}|_2 |\nabla\overline{c}^{k+1}|_2,\\
%\end{split}
%\end{equation}
%\begin{equation}
%\begin{split}\label{litongg3hh}
J_6=&\int_{\mathbb{R}^3} D^\zeta\big(\psi^{k+1}\cdot Q(c^{k+1},\overline{u}^k)\big)\cdot D^\zeta\overline{u}^{k+1}  \text{d}x\\
\leq& C(1+|\overline{E}(c^{k+1})|_\infty)\Big(|\nabla \psi^{k+1}|_2 |\nabla \overline{u}^{k}|_{6} |\nabla \overline{u}^{k+1}|_3+|\psi^{k+1}|_6 |\overline{u}^{k}|_{D^2} |\nabla \overline{u}^{k+1}|_3\Big)\\
&+C|\psi^{k+1}|_6|\nabla c^{k+1}|_6|\nabla \overline{u}^k|_6|\nabla \overline{u}^{k+1}|_2,\\
J_7=&\int_{\mathbb{R}^3} D^\zeta\big(\overline{\psi}^{k+1}\cdot Q(c^k,u^{k-1})\big)\cdot D^\zeta\overline{u}^{k+1} \text{d}x\\
 \leq& C(1+|\overline{E}(c^k)|_\infty)|\overline{\psi}^{k+1}|_2 |\nabla u^{k-1}|_\infty |\nabla D^\zeta\overline{u}^{k+1}|_2,\\
\end{split}
\end{equation}
and
\begin{equation}\label{litongg3}
\begin{split}
J_8=&\int_{\mathbb{R}^3} D^\zeta\Big(\psi^{k+1}(\overline{E}(c^{k+1})-\overline{E}(c^{k}) ) \text{div}u^{k-1}\Big)\cdot D^\zeta\overline{u}^{k+1} \text{d}x\\
\leq & C|\psi^{k+1}|_6 |\nabla^2 \overline{u}^{k+1}|_2|\overline{c}^{k+1}|_3 |\text{div}u^{k-1}|_\infty.
\end{split}
\end{equation}
According to Young's inequality and (\ref{litong3})-(\ref{litongg3}), we have
\begin{equation}\label{gogo12}\begin{split}
&\frac{d}{dt}| \nabla\overline{u}^{k+1}|^2_2+\alpha |\overline{u}^{k+1} |^2_{D^2}\\
\leq& F^k_\eta(t)|\nabla\overline{u}^{k+1}|^2_2+F^k_2(t)\|\overline{c}^{k+1}\|^2_{1}
+F^k_3(t)|\overline{\psi}^{k+1}|^2_{2}+\eta\|\nabla\overline{u}^{k}\|^2_1,
\end{split}
\end{equation}
where
\begin{equation*}
\begin{cases}
F^k_\eta(t)=C\Big(1+\|\nabla c^{k+1}\|^4_1+\frac{1}{\eta^2}(1+\|u^{k}\|^4_{2}+\| u^{k-1} \|^4_{2}+|\psi^{k+1}|^4_{D^1}+|\psi^{k+1}|^2_6|\nabla c^{k+1}|^2_6)\Big), \\[10pt] 
F^k_2(t)=C\Big(\|c^{k+1}\|_2+\|c^{k}\|_2+\|u^k\|_3+|\psi^{k+1}|_6|\text{div} u^{k-1}|_\infty\Big)^2,\ \ 
F^k_3(t)=C\|\nabla u^{k-1}\|^2_2 ,
\end{cases}
\end{equation*}
and we  have  $\int_0^t \big(F^k_\eta(s)+F^k_2(s)+F^k_3(s)\big)\text{d}s\leq \widehat{C}+\widehat{C}_\eta t$ for $t\in (0,T_\epsilon]$.

Then combining (\ref{gogo1}) and (\ref{gogo12}), we easily have
\begin{equation}\label{gogo13}\begin{split}
&\frac{d}{dt}\|\overline{u}^{k+1}\|^2_1+\alpha \|\nabla \overline{u}^{k+1} \|^2_{1}\\
\leq& \Theta^k_\eta(t)\|\overline{u}^{k+1}\|^2_1+\Theta^k_2(t)\|\overline{c}^{k+1}\|^2_{1}
+\Theta^k_3(t)|\overline{\psi}^{k+1}|^2_{2}+\eta\|\nabla\overline{u}^{k}\|^2_1,
\end{split}
\end{equation}
and we  also have  $\int_0^t \big(\Theta^k_\eta(s)+\Theta^k_2(s)+\Theta^k_3(s)\big)\text{d}s\leq \widehat{C}+\widehat{C}_\eta t$, for $t\in (0,T_\epsilon]$.

Finally, let
\begin{equation*}\begin{split}
\Gamma^{k+1}=&\|\overline{c}^{k+1}\|^2_{1}+|\overline{\psi}^{k+1}|^2_{ 2}+\|\overline{u}^{k+1}\|^2_1,
\end{split}
\end{equation*}
then we have
\begin{equation*}\begin{split}
&\frac{d}{dt}\Gamma^{k+1}+\mu\|\nabla\overline{u}^{k+1} \|^2_1
\leq \Pi^k_\eta \Gamma^{k+1}+C\eta \|\nabla \overline{u}^k\|^2_1,
\end{split}
\end{equation*}
for some $\Pi^k_\eta$ such that  $\int_{0}^{t}\Pi^k_\eta(s)\text{d}s\leq \widehat{C}+\widehat{C}_\eta t$.
According to Gronwall's inequality, we have
\begin{equation*}\begin{split}
&\Gamma^{k+1}+\int_{0}^{t}\mu\|\nabla\overline{u}^{k+1}\|^2_1\text{d}s
\leq  \Big( C\eta\int_{0}^{t}  \|\nabla \overline{u}^k\|^2_1\text{d}s\Big)\exp{(\widehat{C}+\widehat{C}_\eta t)}.
\end{split}
\end{equation*}
We can choose $\eta>0$ and $\dot{T}\in (0,T^\epsilon)$ small enough such that
$$
C\eta\exp{\widehat{C}}=\frac{\mu}{4}, \quad \text{and}\quad \text{exp}(\widehat{C}_\eta \dot{T})=2.
$$
Then we easily have
\begin{equation*}\begin{split}
\sum_{k=1}^{\infty}\Big( \sup_{0\leq t\leq \dot{T}} \Gamma^{k+1}+\int_{0}^{\dot{T}} \mu\|\nabla\overline{u}^{k+1}\|^2_1\text{d}s\Big)\leq \widehat{C}<+\infty,
\end{split}
\end{equation*}
%Due to
%\begin{equation*}\begin{split}
%\lim_{k\mapsto \infty} |\overline{\psi}^{k+1}|_{6}\leq C\lim_{k\mapsto \infty}\big( |\overline{\psi}^{k+1}|^{\frac{1}{2}}_{D^1} |\overline{\psi}^{k+1}|^{\frac{1}{2}}_{2}\big)\leq C\lim_{k\mapsto \infty} |\overline{\psi}^{k+1}|^{\frac{1}{3}}_{2}=0,
%\end{split}
%\end{equation*}
which means  that the full consequence $(c^k,\psi^k,u^k)$ converges to a limit $(c,\psi,u)$ in the following strong sense:
\begin{equation}\label{str}
\begin{split}
&c^k\rightarrow c\ \text{in}\ L^\infty([0,\dot{T}];H^1(\mathbb{R}^3)),\\
&\psi^k\rightarrow\psi \ \text{in}\ L^\infty([0,\dot{T}];L^2(B_R)),\\
&u^k\rightarrow u\ \text{in}\ L^\infty ([0,\dot{T}];H^1(\mathbb{R}^3)) \cap L^2([0,\dot{T}];D^2(\mathbb{R}^3)),
\end{split}
\end{equation}
where $B_R$ is a ball centered at origin  with radius $R$, and $R>0$ can be arbitrarily large.

Due to the local estimate (\ref{jkk}) and the lower-continuity of norm for weak or weak$^*$ convergence, we also have $(c,\psi,u)$ satisfies the estimate (\ref{jkk}).
%\begin{equation}\label{gok}\begin{split}
%&\sup_{t\in[0,T]}\Big(\|\rho(t)\|_{H^1\cap W^{1,q}}+ \|\rho_t(t)\|_{L^{2}\cap L^{q}}\Big)+\sup_{t\in[0,T]}\Big(\|(\theta(t),u(t))|_{D^1\cap D^2}\Big)\\
%&+\|I\|_{L^2(\mathbb{R}^+\times S^2;C([0,T_*];H^1\cap W^{1,q}))}+ \|I_t\|_{L^2(\mathbb{R}^+\times S^2;C([0,T_*];L^{2}\cap L^{q}))}\\
%&+\int_{0}^{T}\Big(\|(\theta,u)\|^2_{D^{2,q}}+\|(\theta_t,u_t)\|^2_{D^1}\Big)\text{d}s+\text{ess}\sup_{t\in[0,T]}\Big|\Big(\sqrt{\rho}\theta(t),\sqrt{\rho}u(t)\Big)\Big|_2 \leq C
%\end{split}
%\end{equation}
According to the strong convergence in (\ref{str}), it is easy to see that $(c,\psi,u)$ is a weak solution in the  distribution sense with the regularity (\ref{wode}).
%\begin{equation}\label{rjkqq}\begin{split}
%& c \in L^\infty([0,\dot{T}];H^2),\quad  c_t \in L^\infty([0,\dot{T}];H^1),\ \psi \in L^\infty([0,\dot{T}] ; D^1),\\
%&\partial_i \psi^{(j)}=\partial_j \psi^{(i)} \ (i,j=1,2,3), \ \psi_t \in L^\infty([0,\dot{T}]; L^2),\\
%& u\in L^\infty([0,\dot{T}]; H^2)\cap L^2([0,\dot{T}] ; H^3), \ u_t \in L^\infty([0,\dot{T}]; L^2)\cap L^2([0,\dot{T}] ; D^1),\\
%& t^{\frac{1}{2}}u_{tt}\in L^2([0,\dot{T}];L^2),\quad  t^{\frac{1}{2}}u\in L^\infty([0,\dot{T}];D^3),\\
%&t^{\frac{1}{2}}u_t\in L^\infty([0,\dot{T}];D^1)\cap L^2([0,\dot{T}];D^2),\ t^{\frac{1}{2}}u_{tt}\in L^2([0,\dot{T}];L^2).
%\end{split}
%\end{equation}
So we have given the existence of the strong solution.

\underline{Step 2}. Uniqueness.   Let $(c_1,\psi_1,u_1)$ and $(c_2,\psi_2,u_2)$ be two strong solutions to  Cauchy problem (\ref{li4})-(\ref{vg})  satisfying the uniform a prior estimate (\ref{jkk}). We denote that
$$
\overline{c}=c_1-c_2,\quad \overline{\psi}=\psi_1-\psi_2,\quad \overline{u}=u_1-u_2.
$$
Then according to (\ref{eq:cccq}), $(\overline{c},\overline{\psi},\overline{u})$ satisfies the following system
 \begin{equation}
\label{zhuzhu}
\begin{cases}
\ \overline{c}_t+u_1\cdot \nabla\overline{c} +\overline{u}\cdot\nabla c_{2}+\frac{\gamma-1}{2}(\overline{c} \text{div}u_2 +c_{1}\text{div}\overline{u})=0,\\[10pt]
\ \overline{\psi}_t+\sum_{l=1}^3 A_l(u^1) \partial_l\overline{\psi}+B(u^1)\overline{\psi}+\nabla \text{div}\overline{u}^k=\overline{\Upsilon}_1+\overline{\Upsilon}_2,\\[10pt]
\ \overline{u}_t+ u_1\cdot\nabla \overline{u}+ \overline{u}\cdot \nabla u_{2}+\theta\nabla ((c_1)^2-(c_2)^2) \\[8pt]
=-L(c_1)\overline{u}+\text{div}((\overline{E}(c_1)-(\overline{E}(c_2))\text{div}u_2 \mathbb{I}_3)\\[10pt]
\ \ +\psi_1\cdot Q(c_1,\overline{u})+\overline{\psi}\cdot Q(c_2,u_2)+\psi_1(\overline{E}(c_1)-\overline{E}(c_2))\text{div}u_2,
\end{cases}
\end{equation}
where $\overline{\Upsilon}_1$  and $\overline{\Upsilon}_2$ are defined via
\begin{equation*}
\overline{\Upsilon}_1=-\sum_{l=1}^3(A_l(u^1) \partial_l\psi^{2}-A_l(u^{2}) \partial_l\psi^{2}),\quad \overline{\Upsilon}_2=-(B(u^1) \psi^{2}-B(u^{2}) \psi^{2}).
\end{equation*}

Via the same method used in the derivation of (\ref{go64})-(\ref{gogo1}), let
$$
\Phi(t)=\|\overline{c}(t)\|^2_{1}+|\overline{\psi}(t)|^2_{ 2}+\|\overline{u}(t)\|^2_1,
$$
we similarly have
\begin{equation}\label{gonm}\begin{cases}
\displaystyle
\frac{d}{dt}\Phi(t)+C\|\nabla \overline{u}(t)\|^2_1\leq G(t)\Phi (t),\\[10pt]
\displaystyle
\int_{0}^{t}G(s)ds\leq \widehat{C}\quad  \text{for} \quad 0\leq t\leq \dot{T}.
\end{cases}
\end{equation}
 Then via Gronwall's inequality, the uniqueness follows from $\overline{c}=\overline{\psi}=\overline{u}=0$.

\underline{Step 3}. The time-continuity of the classical solution. It can be obtained via the standard method used in the proof of Lemma \ref{lem1} (see \cite{CK}).

\end{proof}

\subsection{Proof of Remark \ref{r2}}

In this subsection, we  will make a brief disscussion on  the case  $\lambda(\rho)=\rho^b$ when $b\in (1,2) \cup (2,3)$. Here $E(\rho)=\rho^{b-1}$ does not belong to $C^2(\overline{\mathbb{R}}^+)$.

Similarly to the case shown in Theorem \ref{th2}, via introducing  new variables $c$, $\psi$ and  $E(\rho)=\rho^{b-1}$, we need to consider the following Cauchy problem:

\begin{equation}
\begin{cases}
\label{eqedc}
\displaystyle
c_t+u\cdot \nabla c+\frac{\gamma-1}{2}c\text{div} u=0,\\[8pt]
\displaystyle
E_t+u\cdot E+(b-1)E \text{div} u=0,\\[8pt]
\displaystyle
u_t+u\cdot\nabla u +\frac{2}{\gamma-1}c\nabla c+Lu=\psi\cdot Q(c,u),\\[8pt]
(c,E,u)|_{t=0}=(c_0,E_0,u_0),\quad x\in \mathbb{R}^3,\\[8pt]
(c,E,u)\rightarrow (0,0,0) \quad \text{as } \quad |x|\rightarrow \infty,\quad t> 0.
 \end{cases}
\end{equation}
The corresponding existence conclusion can be given as:
\begin{theorem}[\textbf{Existence of the unique local regular solution}]\label{th00}\ \\[2pt]
Let  $1<\gamma \leq 3$. If the initial data $( c_0, E_0, u_0)$ satisfies the regularity conditions
\begin{equation}\label{th00sd}
\begin{split}
&c_0\geq 0,\quad (c_0,E_0, u_0)\in H^2, \quad \psi_0\in  D^1,
\end{split}
\end{equation}
 then there exists a  time $T_*>0$ and a unique regular solution $(c, E, u)$ to Cauchy problem (\ref{eq:1.1})-(\ref{eq:2.211}) with additional regularities:
$$
E \geq 0,\quad  E\in C([0,T_*]; H^2), \quad  E_t \in C([0,T_*]; H^1).
$$
Moreover,  we  have $\rho(t,x)\in C([0,T_*]\times \mathbb{R}^3)$.
\end{theorem}

\begin{proof}
According to the proof  of Theorem \ref{th2} in Subsections $3.1$-$3.4$, the assumptions 
$$
E(\rho)\in C^2(\overline{\mathbb{R}}^+),\quad \text{and} \quad 1< \gamma \leq 2, \quad \text{or} \quad \gamma=3
$$
are only used to deduce the following estimates (see (\ref{jkk})):
$$
|\overline{E}(c)(t)|^2_\infty+\|\overline{E}(c)(t)-\overline{E}(c^\infty)\|^2_{2}+\|\overline{E}(c)_t(t)\|^2_{1}\leq M(b_0)b^4_3,
$$
in Subection $3.2$, and 
$$
\|\overline{E}(c^{k+1})-\overline{E}(c^{k})\|_1\leq C(b_0, \alpha, \gamma,A,T)
$$
in Subsection $3.4$,
where
$$\overline{E}(c)=E(\rho)=E\big(((A\gamma)^{\frac{-1}{2}}c)^{\frac{2}{\gamma-1}}\big) \in C^2(\overline{\mathbb{R}}^+).$$

Thus the key point of our proof for this theorem is to make sure that the desired estimates as above for $E=\rho^{b-1}$ is still avalable based on the additional assumption $E_0\in H^2$.

However, because equations $(\ref{eqedc})_1$ and $(\ref{eqedc})_2$ have totally the same mathematical structure (scalar transport equation),  the desired estimates as above  for $E(\rho)$ can be otained 
via the completely same arguments  used for $c$ as in Subsections $3.1$-$3.4$.

Based on this observation,  we can prove this theorem via the similar arguments used in the proof of Theorem \ref{th2}. Here we omit it.
\end{proof}

\section{Existence of the local strong solution}

Based on the conclusions obtained on Theorem \ref{th2}, we will give the proof for the local existence of strong solutions to the original Cauchy problem (\ref{eq:1.1})-(\ref{eq:2.211}).

\begin{proof} We first give the proof for the case $ 1< \gamma \leq 2 $.
From  Theorem \ref{th2}, we know
there exists  a time $T_* > 0$ such that the Cauchy problem  has a unique regular solution $(c,\psi,u)$ satisfying the regularity (\ref{reg11}), which means that
\begin{equation}\label{reg2}
\begin{split}
(\sqrt{A\gamma}\rho^{\frac{\gamma-1}{2}},u )=(c,u)\in C((0,T_*)\times \mathbb{R}^3).
\end{split}
\end{equation}
 According to  transformation
$$
\rho(t,x)=\Big(\frac{c}{\sqrt{A\gamma}}\Big)^{2\theta}(t,x),
$$
and $ 2\theta\geq 2$ due to $ 1< \gamma \leq 2 $,  it is easy to show that
$$\rho(t,x)\in C((0,T_*)\times\mathbb{R}^3)\cap C([0,T_*];H^2).$$

Multiplying $(\ref{eq:cccq})_1$ by
$
\frac{\partial \rho}{\partial c}(t,x)=\frac{2\theta}{\sqrt{A\gamma}} \Big(\frac{c}{\sqrt{A\gamma}}\Big)^{2\theta-1}(t,x)\in C((0,T_*)\times \mathbb{R}^3)
$,
we get the continuity equation $(\ref{eq:1.1})_1$:
\begin{equation} \label{eq:2.58}
\rho_t+u \cdot\nabla \rho+\rho\text{div} u=0.
\end{equation}
Then combining (\ref{eq:2.58}) and
$
u(t,x)\in C([0,T_*],H^{2})\bigcap C^1([0,T_*],H^{1})$, from the linear qusi-linear hyperbolic equation theory, we  immediately have
$$
\rho \in C([0,T_*],H^{2})\cap C^1([0,T_*],H^{1}).
$$

Multiplying $(\ref{eq:cccq})_2$ by
$
\Big(\frac{c}{\sqrt{A\gamma}}\Big)^{2\theta}=\rho(t,x)\in C((0,T_*)\times \mathbb{R}^3)
$,
we get the momentum equations $(\ref{eq:1.1})_2$:
\begin{equation} \label{eq:2.60}
\begin{split}
&\rho u_t+\rho u\cdot \nabla u+\nabla P=\text{div}\Big(\alpha\rho(\nabla u+ (\nabla u)^\top)+\rho E(\rho)\text{div}u I_3\Big).
\end{split}
\end{equation}
That is to say, $(\rho,u)$ satisfies the compressible isentropic Navier-Stokes equations (\ref{eq:1.1})   a.e. in  $ (0,T_*]\times \mathbb{R}^3$ and has the regularity (\ref{reg11}) with $$
\rho \in C([0,T_*],H^{2})\cap C^1([0,T_*],H^{1}).$$

From the continuity equation and Lemma $6$ in \cite{CK}, it is easy to get that
the solution $\rho$ is represented by the formula
$$
\rho(t,x)=\rho_0(U(0;t,x))\exp\Big(\int_{0}^{t}\textrm{div} u(s,U(s;t,x))\text{d}s\Big),
$$
which, together with  $\rho_0\geq 0$, immediately implies that 
$$
\rho(t,x)\geq 0, \ \forall (t,x)\ \in  \ [0,T_*]\times \mathbb{R}^3.
$$
In summary, the Cauchy problem (\ref{eq:1.1})-(\ref{eq:2.211}) has a unique strong solution $(\rho,u)$.

Finally, when $\gamma=3 $, we quickly have the relation  $\rho(t,x)=\frac{1}{\sqrt{A\gamma}}c(t,x)$, via the same argument used in the case  $1< \gamma \leq 2$ as above, the same conclusions will be obtained.

\end{proof}

%Actually, for any positive integer $k\geq 0$, via the similar arguments, we have the following result:
%\begin{theorem}\label{ab22}
%If the initial data $( \rho_0, u_0)$ satisfies the regularity condition
%\begin{equation}\label{pth78rr}
%\begin{split}
%&\rho_0\geq 0,\quad (\rho^{\frac{\gamma-1}{2}}_0, u_0)\in H^{3+k}, \quad \nabla \rho_0/\rho_0\in  D^1\cap D^{2+k}
%\end{split}
%\end{equation}
%and the compatibility condition
%\begin{equation}\label{th79rr}
%\begin{split}
%2\theta\rho^{\frac{\gamma-1}{2}}_0 \nabla \rho^{\frac{\gamma-1}{2}}_0+Lu_0-(\nabla \rho_0/\rho_0)\cdot Q(u_0)=g_1
%\end{split}
%\end{equation}
%for some $g_1\in H^1$.
%Then there exists a small time $T_*$ and a unique regular solution $(\rho, u)$ to Cauchy problem (\ref{eq:1.1})-(\ref{eq:2.211}) satisfying
%\begin{equation*}\begin{split}
%& \rho^{\frac{\gamma-1}{2}} \in C([0,T_*];H^{3+k}),\quad (\rho^{\frac{\gamma-1}{2}})_t \in C([0,T_*];H^{2+k}) ,\ \nabla \rho/\rho \in C([0,T_*] ; D^1\cap D^{2+k}),\\
%&  u\in C([0,T_*]; H^{3+k})\cap L^2([0,T_*] ; H^{4+k}), \ u_t \in C([0,T_*]; H^{1+k})\cap L^2([0,T_*] ; D^{2+k}).
%\end{split}
%\end{equation*}
%Moreover, if $1< \gamma \leq 3$, we also have $\rho(t,x)\in C^1([0,T_*]\times \mathbb{R}^3)$, which means that the solution $(\rho, u)$ satisfies the equations (\ref{eq:1.1}) classically. And for any positive integer $j\geq 1$, we can also  have
%\begin{equation*}\begin{split}
%& \rho \in C([0,T_*];H^{j}),\quad \rho_t \in C([0,T_*];H^{j-1}),\quad \text{for} \ 1< \gamma \leq 1+\frac{2}{j},\ \text{or}\ \gamma=2.\end{split}
%\end{equation*}
%\end{theorem}

\section{No-existence of global solutions with $L^\infty$ decay on $u$}

In order to prove the phenomenon shown in Theorem \ref{th:2.20}, firstly  we need to  introduce some physical notations:
\begin{align*}
&m(t)=\int_{\mathbb{R}^{3}}\rho(t,x)\text{d}x \quad \textrm{(total mass)},\\
%&\mathbb{F}=\int_{\mathbb{R}^{3}} \rho(t,x)u(t,x)\cdot x \text{d}x \quad \textrm{ (radial component of momentum)},\\
&E_k(t)=\frac{1}{2}\int_{\mathbb{R}^{3}} \rho(t,x)|u(t,x)|^{2}\text{d}x \quad \textrm{ (total kinetic energy)}.
%&E_i(t)=\int_{\mathbb{R}^{3}} \frac{P(t,x)}{\gamma-1}\text{d}x \quad \textrm{ (total potential energy)},\\
%&M(t)=\int_{\mathbb{R}^{3}} \rho(t,x)|x|^{2}\text{d}x \quad \textrm{ (second moment)}.
%&\varepsilon(t)=\int_{\mathbb{R}^{3}} E(t,x) \text{d}x=E_k(t)+E_i(t) \quad  \textrm{(total  energy)},
\end{align*}
%where the total energy density $E(t,x)=\frac{1}{2}\rho u^2+\frac{p}{\gamma-1}$. Moreover, combining the continuity equation and the momentum equations, we obtain the following relation that $E(t,x)$ has to satisfy
%\begin{equation}\label{ener}
%E_t+\text{div}((E+p)u)=u\text{div} \mathbb{T}.
%\end{equation}

Based on the existence theory established in Theorem \ref{th2} and the additional initial conditions in Theorem \ref{th:2.20}, we can show  that there exists a unique
regular solution
$ (\rho,u)(t,x)$ on $[0,T]\times \mathbb{R}^3$ which has finite mass $m(t)$, finite momentum $\mathbb{P}(t)$, finite kinetic energy $E_k(t)$.
Actually, due to $1<\gamma \leq 2$,  we have 
$$
m(t)=\int_{\mathbb{R}^3} \rho \text{d}x \leq C\int_{\mathbb{R}^3} c^{\frac{2}{\gamma-1}} \text{d}x \leq C|c|^2_2<+\infty,
$$
which, together with the regularity shown in Theorem \ref{th2}, implies that 
\begin{equation}\label{finite}
%\displaystyle
%\mathbb{F}(t)=\int_{\mathbb{R}^3} \rho u \text{d}x \leq |u|_\infty\int_{\mathbb{R}^3} \rho  \text{d}x <+\infty,\\[12pt]
\displaystyle
E_k(t)=\int_{\mathbb{R}^3} \frac{1}{2}\rho|u|^2 \text{d}x \leq C|\rho|_\infty|u|^2_2<+\infty.
\end{equation}

Secondly, we give the following lemmas which are the  revised versions for the constant viscosity case \cite{olg1}.
\begin{lemma}
\label{lemmak} Let $1< \gamma \leq 2$ and  $(\rho,u)$ be the regular solution obtained in Theorem \ref{th2} with the 	
additional initial conditions shown in Theorem \ref{th:2.20}, then
$$\mathbb{P}(t)=\mathbb{P}(0), \quad  m(t)= m(0), \quad \text{for} \quad t\in [0,T]. $$
\end{lemma}
\begin{proof}
 According to the momentum equations, we immediately deduce that
\begin{equation}\label{deng1}
\mathbb{P}_t=-\int_{\mathbb{R}^3} \text{div}(\rho u \otimes u)\text{d}x-\int_{\mathbb{R}^3} \nabla P\text{d}x+\int_{\mathbb{R}^3} \text{div}\mathbb{T}\text{d}x.
\end{equation}

We first claim that   
$$\int_{\mathbb{R}^3} \text{div}\mathbb{T}\text{d}x=0.$$
Let $R> 0$ be a arbitrarily large constant, from Green's formula, we only need to prove
\begin{equation}\label{deng}
\lim_{R\rightarrow +\infty}\int_{\partial B_R} \mathbb{T}\cdot n \text{dS}=\lim_{R\rightarrow +\infty}\int_{\partial B_R} \rho(\alpha(\nabla u+(\nabla u)^\top)+E(\rho)\text{div}u\mathbb{I}_3)\cdot n \text{dS}=0.
\end{equation}
We denote 
$$G_R=\Big|\int_{\partial B_R} \rho \nabla u\cdot n \text{dS}\Big|.$$ According to Definition \ref{d1},  we have
$$\rho \in C([0,T]; H^2),\quad \nabla u\in C([0,T]; H^1),$$ from H\"older's inequality, which implies that
\begin{equation}\label{deng2}
\int_{\mathbb{R}^3}\rho |\nabla u| \text{d}x\leq |\rho|_2 |\nabla u|_2< \infty, \quad \text{for}\ t\in [0,T].
\end{equation}

Next let $\Omega_1=B_1$, $\Omega_i=B_i/B_{i-1}$ $(i\geq 2)$, from (\ref{deng2}), we have
\begin{equation}\label{deng3}
\int_{\mathbb{R}^3}\rho |\nabla u| \text{d}x=\sum_{i=1}^{\infty} \int_{\Omega_i}\rho |\nabla u| \text{d}x< \infty, \quad \text{for}\ t\in [0,T].
\end{equation}
Then we immediately obtain that
\begin{equation}\label{deng4}
\lim_{i\rightarrow \infty}\int_{i-1}^{i} G_R \text{d}R \leq \lim_{i\rightarrow \infty}\int_{\Omega_i}\rho |\nabla u| \text{d}x=0.
\end{equation}

Next we prove that $G_R$ is a uniformly continuous function with respect to $R$,  let $0< R_1< R_2<\infty$ be two constants, we have
\begin{equation}\label{deng5}
\begin{split}
|G_{R_1}-G_{R_2}|\leq& \Big|\int_{\partial (B_{R_2}/B_{R_1})} \rho \nabla u \cdot n \text{dS} \Big|\\
=& \Big|\int_{B_{R_2}/B_{R_1}} \text{div}(\rho \nabla u ) \text{dx} \Big|
 \leq  \|\rho\|_{W^{1,6}} \|\nabla u\|_1 |B_{R_2}/B_{R_1}|^{\frac{1}{3}},
\end{split}
\end{equation}
where $|B_{R_2}/B_{R_1}|$ is the three-dimensional Lebesgue measure.

 At last, if 
$$\lim_{R\rightarrow +\infty} G_R\neq 0,$$ we know that there exists a constant $\epsilon_0 > 0$, for arbitrarily large $R> 0$, there exists a constant $R_0 > R$ such that $G_{R_0} \geq \epsilon_0 $. Due to the uniform continuity, we know that there exists a small constant $\eta > 0$ such that
$$ |G_{R_0}-G_R|\leq \frac{\epsilon_0}{2}  \quad \text{for} \quad  |R_0-R|\leq \eta,$$
which means that
\begin{equation}\label{deng6}
G_R\geq \frac{\epsilon_0}{2},  \quad \text{for} \quad  |R_0-R|\leq \eta.
\end{equation}
It is obvious that, for  sufficiently large $i$, there always exists some  $j\geq i$ such that
\begin{equation}\label{deng7}\int_{j-1}^{j} G_R \text{d}R\geq\frac{\eta\epsilon_0}{2},\end{equation}
which is impossible due to (\ref{deng4}). So we immediately have that 
$$\lim_{R\rightarrow +\infty} G_R=0,$$
 which makes sure  (\ref{deng}) holds. Then via the similar arguments used to prove  (\ref{deng}), we also can deduce that 
$$
-\int_{\mathbb{R}^3} \text{div}(\rho u \otimes u)\text{d}x-\int_{\mathbb{R}^3} \nabla P\text{d}x=0,
$$
which, together with (\ref{deng1})-(\ref{deng}), immediately implies the conservation the momentum.

Similarly, we also can get the conservation of mass, the  proof is similar without essential modifications, here we omit it.
\end{proof}
\begin{lemma}\label{mo}
Let $1< \gamma \leq 2$ and  $(\rho,u)$ be the regular solution obtained in Theorem \ref{th2} with the 	
additional initial conditions shown in Theorem \ref{th:2.20}, there exists a unique lower bound $C_0$ which has no dependent on $t$ for $E_k(t)$ such that
$$
E_k(t)\geq C_0> 0 \quad \text{for} \quad t\in [0,T].
$$
\end{lemma}
\begin{proof}
Due to H\"older's inequality and momentum equations, we deduce that
\begin{equation}\label{rty}
\begin{split}
|\mathbb{P}(0)|=&|\mathbb{P}(t)|\leq\int_{\mathbb{R}^3}\rho(t,x) |u|(t,x) \text{d}x\\
  \leq& \sqrt{2}m^{\frac{1}{2}}(t)E^{\frac{1}{2}}_k(t)=\sqrt{2}m^{\frac{1}{2}}(0)E^{\frac{1}{2}}_k(t),
\end{split}
\end{equation}
which implies that there exists a unique positive lower bound for $E_k(t)$ such that
\begin{equation}\label{rty1}
E_k(t)\geq \frac{|\mathbb{P}(0)|^2}{2m(0)}> 0  \quad \text{for} \quad t\in [0,T].
\end{equation}
\end{proof}
\begin{remark}\label{gu66}
The positive lower bound of the total kinetic energy $E_k(t)$ will play an key role in the proof of the corresponding non-existence of global regular solutions with $L^\infty$ decay on $u$, which is essentially obtained via the conservation of the momentum based on the regularity of regular solutions. The same conclusions can't be obtained for the strong solutions shown in  \cite{CK3} or \cite{CK} because of the  different mathematical structure,  even if the initial mass density and velocity are both compactly supported. In this sense, the definition of regular solutions with vacuum is consistent with the physical background of the compressible Navier-Stokes equations.
\end{remark}
%\begin{proof}
%According to  Cauchy's inequality and the definition of the stress tensor $\mathbb{T}$, we easily know that
%\begin{equation}\label{eq:5.9}
%\begin{split}
%\nabla \cdot (u \mathbb{T })
%=&\sum_{i=1}^{3}\sum_{j=1}^{3}\left( u_{i}\frac{\partial{\mathbb{T }_{ij}}}{\partial{x_j}}+\partial_{j}u_{i}\mathbb{T }_{ij}\right)\\
%=&u\cdot (\nabla \cdot \mathbb{T }) +\sum_{i=1}^{3}\sum_{i=1}^{3} \partial_{j}u_{i} \Big(\mu(\rho)  (\partial_{j}u_{i}+\partial_{i}u_{j})+\lambda(\rho)\delta_{ij}\nabla \cdot u\Big)\\
%=&u\cdot (\nabla \cdot \mathbb{T })+2\mu(\rho)\sum_{i=1}^{3}(\partial_{i}u_{i})^{2}+\mu(\rho)  \sum_{i\neq j}^{3} (\partial_{i}u_{j})^2\\
%&+2\mu(\rho) \sum_{i>j} (\partial_{i}u_{j})(\partial_{j}u_{i})+\lambda(\rho) \left(\sum_{i=1}^{3}\partial_{i}u_{i}\right)^{2}\\
%\geq& u\cdot (\nabla \cdot \mathbb{T })+\big(\lambda(\rho)+\frac{2}{3}\mu(\rho)\big)(\nabla\cdot u)^2
%=u\cdot (\nabla \cdot \mathbb{T })+\big(\alpha+\frac{2}{3}\beta)\rho (\nabla\cdot u)^2.
%\end{split}
%\end{equation}
%Combining the regularity of $u$, (\ref{ener}) and (\ref{eq:5.9}), we easily obtain the desired conclusion.
%\end{proof}

\textbf{Next we give the proof for Theorem \ref{th:2.20}:}
\begin{proof}
Combining the definition of $E_k(t)$ and Lemmas \ref{lemmak}-\ref{mo}, we easily have
$$
C_0\leq E_k(t)\leq \frac{1}{2} m(0)|u(t)|^2_\infty \quad \text{for} \quad t\in [0,T],
$$
which means that there exists a positive constant $C_u$ such that
$$
|u(t)|_\infty\geq C_u  \quad \text{for} \quad t\in [0,T].
$$
Then we quickly obtain the desired conclusion as shown in Theorem \ref{th:2.20}.
\end{proof}

{\bf Acknowledgement:} The research of  S. Zhu was supported in part
by National Natural Science Foundation of China under grant 11231006, Natural Science Foundation of Shanghai under grant 14ZR1423100 and  China Scholarship Council.


\bigskip

\begin{thebibliography}{99}




 \bibitem{bd2} D. Bresch,  B. Desjardins and C. Lin,  On some compressible fluid models: Korteweg, Lubrication, and Shallow water systems, \textit{Commun. Part. Differ. Equations}, \textbf{28} (2003), 843-868.
 \bibitem{bd} D. Bresch,  B. Desjardins and G. M$\acute{e}$tivier,  Recent mathematical results and open problems about shallow water equations, \textit{Anal. Simu. Fluid. Dynam}, (2006), 15-31.

 \bibitem{bd3} D. Bresch  and  B. Desjardins,  Some diffusive capillary  models of Korteweg type, \textit{C. R.  Acad. Science}, Vol. \textbf{332} No. \textbf{11} (2004), 881-886.









\bibitem{CK3} Y. Cho, H.  Choe  and H. Kim, Unique solvability of the initial boundary value problems for compressible viscous  fluids, \textit{J. Math. Pures  Appl.}, \textbf{83} (2004), 243-275.



\bibitem{CK} Y. Cho and  H. Kim, Existence results for  viscous polytropic fluids with vacuum, \textit{J. Differential Equations}, \textbf{228} (2006), 377-411.

\bibitem{fu1} E. Feireisl, A. Novotny and  H. Petzeltová, On the existence of globally defined weak solutions to the Navier-Stokes equations, J\text{. Math. Fluid Mech.},  \textbf{3(4)}  (2001), 358-392.

\bibitem{fu2} E.  Feireisl,  On the motion of a viscous, compressible, and heat conducting fluid, \textit{ Indiana Univ. Math. J.},  \textbf{ 53(6)} (2004), 1705-1738.

\bibitem{fu3} Feireisl,  \textit{Dynamics of Viscous Compressible Fluids}, Oxford: Oxford University Press, 2004.




\bibitem{gandi} G.  Galdi, \textit{An introduction to the Mathmatical Theory of the Navier-Stokes equations}, Springer, New York, 1994.

\bibitem{hoff} D. Hoff and  D. Serre,  The failure of continuous dependence on initial data for the   Navier-Stokes equations for compressible flow, \textit{ SIAM J. Appl. Math.}, \textbf{51} (1991), 887-898.




\bibitem{HX1} X. Huang, J. Li and  Z. Xin, Global Well-posedness of classical solutions with large oscillations and vacuum to the Three-Dimensional Isentropic Compressible Navier-Stokes Equations, \textit{ Comm. Pure  Appl. Math.}, \textbf{65} (2012), 549-585.







\bibitem{sz3} Y. Li, R. Pan and S. Zhu, 2D compressible Navier-Stokes equations with degenerate viscosities and far field vacuum, (2013) Submitted.

        \bibitem{sz345} Y. Li, R. Pan and S. Zhu,  On regular solutions for  viscous polytropic fluids  with degenerate viscosities  and vacuum, (2014) Preprint. 




  \bibitem{sz2} Y. Li and S. Zhu, Formation of singularities in solutions to the compressible radiation hydrodynamics equations  with vacuum,  \textit{J. Differential Equations},  \textbf{256} (2014), 3943-3980.
 \bibitem{dcds} Y. Li and S. Zhu, On regular solutions of the 3-D compressible isentropic Euler-Boltzmann equations with vacuum, to appear in \textit{Discrete Contin. Dynam. Systems},  (2014),  Accepted.  


  \bibitem{oar} O.  Ladyzenskaja and N.  Ural'ceva,\textit{ Linear and Quasilinear Equations of Parabolic Type}, American Mathematical Society, Providence, RI, 1968.
  \bibitem{hailiang} H. Li, J. Li and Z. Xin, Vanishing of vacuum states and blow-up phenomena of the compressible Navier-Stokes equations, \textit{Commun. Math. Phys.}, \textbf{281} (2008), 401-444.

\bibitem{tlt} Tatsien Li  and T. Qin, \textit{Physics and Partial Differential Equations},  Siam: Philadelphia,   Higher Education Press: Beijing, 2014.


  \bibitem{lions} P.  Lions, \textit{Mathematical topics in fluid dynamics} In: Compressible Models. Oxford University Press, \textbf{2} 1998.





\bibitem{tpy} T. Liu and  T. Yang, Compressible Euler equations with vacuum, \textit{J. Differential Equations}, \textbf{140} (1997), 223-237.
\bibitem{zhangjianwen} S, Liu, J. Zhang and  J. Zhao, Global classical solutions for 3D compressible Navier-Stoles equations with vacuum and a density-dependent viscosity coefficient, \textit{ Jour. Math. Anal. Appl.}, \textbf{401} (2013), 795-810.
\bibitem{taiping} T. Liu, Z. Xin and T. Yang, Vacuum states for compressible flow, \textit{Discrete Contin. Dynam. Systems}, \textbf{4} (1998), 1-32.

   \bibitem{luoluo} Z. Luo, Local existence of classical solutions to the two-dimensional viscous compressible flows with vacuum, \textit{Commun. Math. Sci.}, \textbf{10} (2012), 527-554.

\bibitem{amj} A. Majda, \textit{Compressible fluid flow and systems of conservation laws in several space variables}, Applied Mathematical Science \textbf{53}, Spinger-Verlag: New York, Berlin Heidelberg, 1986.


\bibitem{tms1} T. Makino,  S. Ukai and  S. Kawashima, Sur la solution $\grave{\text{a}}$ support compact de equations d'Euler compressible, \textit{Japan J.  Appl.  Math.}, \textbf{33} (1986),  249-257.

\bibitem{makio} T. Makino, On a local existence theorem for the evolution  equation of gaseous stars,\textit{ Transport Theory Statist. Phys.}, \textbf{21} (1992), 615-624.

\bibitem{vassu} A. Mellet and A. Vasseur, On the barotropic compressible Navier-Stokes equations,
\textit{Commun. Part. Differ. Equations},  \textbf{32} (2007), no. 1-3, 431--452.


\bibitem {nash}
  J. Nash,
Le probleme de Cauchy pour les {\'e}quations diff{\'e}rentielles d\'un fluide g{\'e}n{\'e}ral,
  \textit{Bull. Soc. Math. France}, \textbf {90} (1962), 487-491.

\bibitem{olg1} O. Rozanova, Blow-up of smooth highly decreasing at infinity solutions to  the compressible Navier-Stokes Equations, \textit{J. Differential Equations}, \textbf{245} (2008), 1762-1774.




\bibitem{jm} J. Simon, Compact sets in $L^P(0,T;B)$, \textit{Ann. Mat. Pura. Appl.}, \textbf{ 146} (1987), 65-96.


 \bibitem{harmo} E.  Stein, \textit{Singular integrals and Differentiablility properties of Functions}, Princeton Univ. Press, Princeton NJ,  1970.



\bibitem{zx} Z.  Xin, Blow-up of smooth solutions to the compressible Navier-Stokes Equations with Compact Density, \textit{Commun. Pure  Appl. Math.}, \textbf{ 51} (1998), 0229-0240.

\bibitem{zwy} Z.  Xin and W. Yan, On blow-up of classical solutions to the compressible Navier-Stokes Equations, \textit{Comm. Math. Phys.},  \textbf{321} (2013), 529-541.


\bibitem{tyc2} T. Yang and  C. Zhu, Compressible Navier-Stokes equations with degnerate viscosity coefficient and vacuum, \textit{Commun. Math. Phys.}, \textbf{230} (2002), 329-363. 


\bibitem{zyj} T. Yang and  H. Zhao, A vacuum problem for the one-dimensional compressible Navier-Stokes equations with density-dependent viscisity,  \textit{J. Differential Equations}, \textbf{184} (2002), 163-184.


\end{thebibliography}
\end{document}